\documentclass[11pt,
]{article}

\usepackage{amssymb,amsfonts,amstext,amsmath,amsthm,latexsym,mathrsfs}
\usepackage[bbgreekl]{mathbbol}

\usepackage{makeidx,mparhack,vector}

\makeatletter
\if@twoside%
\else\fi%
\makeatother  
\newcommand{\nek}{\newcommand}
\nek{\renek}{\renewcommand}
\nek{\vyk} [1] {}

\nek{\ubf}{\fontseries{b}\selectfont}
\nek{\bfit}{\bfseries\itshape}
\nek{\bftt}{\ttfamily\bfseries\upshape\selectfont}

\nek{\parf}{\subsection}
\renek{\thesubsection}{\arabic{subsection}}

\theoremstyle{plain}

\newtheorem{theorem}             {Theorem} [subsection]
\newtheorem{corollary}  [theorem]{Corollary}
\newtheorem{prop}  [theorem]{Proposition}
\newtheorem{lemma}      [theorem]{Lemma}
\newtheorem{claim}      [theorem]{Claim}

\theoremstyle{definition}

\newtheorem{definition} [theorem]{Definition}
\newtheorem{vopr}       [theorem]{Question}
\newtheorem{prim}       [theorem]{Example}
\newtheorem{zam}        [theorem]{Remark}
\newtheorem*{prF}   {Proof}    
\newtheorem*{ack}   {Acknowledgement}    
\newtheorem{ggi}   [theorem] {Blanket Assumption}  

\nek{\thsp}{\hspace{0.1ex}}
\nek{\bvo}{\begin{vopr}}
\nek{\evo}{\end{vopr}}
\nek{\back}{\begin{ack}}
\nek{\eack}{\end{ack}}
\nek{\bpro}{\begin{prop}}
\nek{\epro}{\end{prop}}
\nek{\bcor}{\begin{corollary}}
\nek{\ecor}{\end{corollary}}
\nek{\bdf} {\begin{definition}}
\nek{\edf} {\qed\end{definition}}
\nek{\eDf} {\end{definition}}
\nek{\bgg} {\begin{ggi}}
\nek{\egg} {\qed\end{ggi}}
\nek{\ble} {\begin{lemma}}
\nek{\ele} {\end{lemma}}
\nek{\bcl} {\begin{claim}}
\nek{\ecl} {\end{claim}}
\nek{\bpri}{\begin{prim}}
\nek{\epri}{\qed\end{prim}}
\nek{\eprI}{\end{prim}}
\nek{\bte} {\begin{theorem}}
\nek{\ete} {\end{theorem}}
\nek{\baq} {\begin{aq}}
\nek{\eaq} {\end{aq}}
\nek{\bre}{\begin{zam}}
\nek{\ere}{\qed\end{zam}}
\nek{\bpf} {\begin{prF}} 
\nek{\epf} {\qed\end{prF}} 
\nek{\epg} {\end{prF}} 
\nek{\epF}[1]{\hfill\hbox{$\square$\ ({\small#1})}\end{prF}}

\nek{\qeD}[1]{\hfill\hbox{$\square$\ ({\small#1})}}

\nek{\ben}{\begin{enumerate}\itemsep=0.2em}
\nek{\een}{\end{enumerate}}
\nek{\bde}{\begin{description}\itemsep=0.2em}
\nek{\ede}{\end{description}}
\nek{\bit}{\begin{itemize}\itemsep=0.2em}
\nek{\eit}{\end{itemize}}
\nek{\bay}{\begin{array}}
\nek{\eay}{\end{array}}
\nek{\bce}{\begin{center}}
\nek{\ece}{\end{center}}
\nek{\bqu}{\begin{quotation}\noindent}
\nek{\equ}{\end{quotation}}
\nek{\lh}  {\mathop{\tt lh}}
\nek{\dom} {\mathop{\tt dom}}
\nek{\clos} {\mathop{\tt clos}}
\nek{\num} {\mathop{\tt num}}
\nek{\ran} {\mathop{\tt ran}}
\nek{\rank} {\mathop{\tt rank}}
\nek{\tsup} {\mathop{\tt sup}}
\nek{\tmax} {\mathop{\tt max}}
\nek{\sep}{{\tt Sep}}
\nek{\red}{{\tt Red}}
\nek{\TC}  {{\rm TC}\hspace{0.4ex}}
\nek{\hc}  {\mathbb{HC}}
\nek{\rL} {{\mathbb L}}
\nek{\rV} {{\mathbb V}}

\nek{\ZFC} {\text{\ubf ZFC}}
\nek{\zfc} {\ZFC}
\nek{\zfcm} {\zfc^-}
\nek{\al}  {\alpha}
\nek{\ga}  {\gamma}
\nek{\Ga}  {\Gamma}
\nek{\da}  {\delta}
\nek{\Da}  {\Delta}
\nek{\kpa} {\kappa}
\nek{\la}  {\lambda}
\nek{\ve}  {\varepsilon}
\nek{\vpi} {\varphi}
\nek{\vpo} {\overline\vpi}
\nek{\sg}  {\sigma}
\nek{\Sg}  {\Sigma}
\nek{\om}  {\omega}
\nek{\Om}  {\Omega}
\nek{\lom} {^{<\om}} 
\nek{\omi} {\om_1}
\nek{\oli} {\omi^\rL}
\nek{\omb} {\om_2}
\nek{\omil}{\om_1^\rL}
\nek{\za}  {\zeta}
\nek{\ali} {\aleph_1}
\nek{\ald} {\aleph_2}
\nek{\alil} {\aleph_1^\rL}
\nek{\fs}[2]{{\mathbf\Sigma}^{#1}_{#2}}
\nek{\fp}[2]{{\mathbf\Pi}^{#1}_{#2}}
\nek{\fd}[2]{{\mathbf\Delta}^{#1}_{#2}}
\nek{\iDa}{\varDelta}
\nek{\iSg}{\varSigma}
\nek{\iPi}{\varPi}
\nek{\is}[2]{\iSg^{#1}_{#2}}
\nek{\ip}[2]{\iPi^{#1}_{#2}}
\nek{\id}[2]{\iDa^{#1}_{#2}}
\nek{\BBB}{\hspace{0.1ex}}
\nek{\dP}{{\BBB{\mathbb P}\BBB}}
\nek{\dQ}{{\BBB{\mathbb Q}\BBB}}
\nek{\gM}{{\BBB{\goth M}\BBB}}
\nek{\gN}{{\BBB{\goth N}\BBB}}
\nek{\cM}{{\BBB{\skri M}\BBB}}
\nek{\cP}{{\BBB{\skri P}\BBB}}
\nek{\pu}  {\varnothing}
\nek{\sq}  {\subseteq}
\nek{\qs}  {\supseteq}
\nek{\su}  {\subset}
\nek{\sneq}{\subsetneqq}
\nek{\eqv} {\mathbin{\,\Longleftrightarrow\,}}
\nek{\imp} {\mathbin{\,\Longrightarrow\,}}
\nek{\mpi} {\mathbin{\,\Longleftarrow\,}}
\nek{\lra} {\longrightarrow} 
\nek{\we}  {{\mathbin{\hspace*{0.2ex}^\wedge}}}
\nek{\sus} {{\exists\,}}
\nek{\kaz} {{\forall\,}}
\nek{\ti}  {\times}
\nek{\dm}  {$$}
\nek{\iy}  {\infty}
\nek{\abs}[1]{|#1|}
\nek{\abt}[1]
{\mathopen{{|}\hspace*{-0.2ex}{|}}#1\mathclose{{|}\hspace*{-0.2ex}{|}}}

\nek{\nin}{\notin}
\nek{\res} {{\hspace*{0.1ex}\restriction\hspace*{0.3ex}}}

\nek{\ang} [1] {\langle #1\rangle}
\nek{\ans} [1] {\{\hspace{0.1ex}#1\hspace{0.1ex}\}}
\nek{\dd}[1]{$\mtho\hspace{0.2ex}{#1}$-\hspace{0.0ex}}
\nek{\itla}{\item\label}
\nek{\mtho}{\mathsurround=0mm}
\nek{\msur}{\hspace*{-1\mathsurround}}
\nek{\dsur}{\hspace{-0.3\mathsurround}}
\nek{\hsur}{\hspace{-0.5\mathsurround}}
\nek{\noi}{\noindent}
\nek{\vom}{\vspace{1mm}}
\nek{\vtm}{\vspace{2mm}}

\nek{\bez}{\smallsetminus}
\nek{\rit} [1] {{\it#1\/}}

\nek{\ens} [2] {\ans{{#1\hspace{0.5ex}{:}}\zz\hspace{0.5ex}#2}}
\nek{\zz} {\linebreak[0]} 
\nek{\goth}{\mathfrak}
\nek{\skri}{\mathscr}

\nek{\yo} {,\linebreak[0]}
\nek{\yi} {\hspace{\mathsurround},\linebreak[0]\hspace{\mathsurround}}
\nek{\yd} {\hspace{\mathsurround},\linebreak[0]\:}
\nek{\yt} {\hspace*{\mathsurround}\text{,}\linebreak[0]\;}

\nek{\viv} {\text{vice versa}}
\nek{\ie} {\text{\sl i.\resp e.}}
\nek{\eg} {\text{\sl e.\resp g.}}
\nek{\poo} {\text{w.\resp r.\resp t.}}
\nek{\noo} {\text{w.\resp l.\resp o.\resp g.}}

\nek{\iesp}{\hspace{0.3ex}}
\nek{\resp}{\hspace{0.25ex}}

\nek{\itsep}{\itemsep=0.25ex plus 0.1ex minus 0.1ex}

\nek{\tenu}[1]{
\def\theenumi{#1}
\itsep}

\nek{\tenv}[1]{
\def\theenumii{#1}
\itsep}

\nek{\Aenu}{\tenu{(\Alph{enumi})}}
\nek{\aenu}{\tenu{(\alph{enumi})}}
\nek{\nenu}{\tenu{{\rm(\arabic{enumi})}}}

\nek{\fenu}{\itsep\tenu{{\mtho$(\fnsymbol{enumi})$}}}
\nek{\renu}{\itsep\tenu{{\rm(\roman{enumi})}}\itsep}
\nek{\Renu}{\itsep\tenu{{\rm(\Roman{enumi})}}\itsep}

\nek{\atc} {\addtocounter{enumi}1}
\nek{\atm} {\addtocounter{enumi}{-1}}

\nek{\stk} [2] {\ang{#1\hspace{0.3ex};\hspace{0.1ex}#2}}
\nek{\sis} [2] {\ang{#1}_{#2}}
\nek{\sid} [3] {\ang{#1}^{#2}_{#3}}

\renek{\cM} {\goth M}

\nek{\lam} [1]
{\label{#1}\hspace*{-3pt}\imar{#1}%
}%
\nek{\las} [1]
{\label{#1}\imae{#1}}%

\nek{\imar}[1]{\marginpar[%\vspace{-1ex}%
\flushright\footnotesize%
$\mtho\longrightarrow$\\%
\vspace{-1ex}{#1}\vspace*{1ex}]%
{%\vspace{-1ex}%
\flushleft\footnotesize%
$\mtho\longleftarrow$\\%
\vspace{-1ex}{#1}\vspace*{1ex}}%
}%
\nek{\imae}[1]{\marginpar[%\vspace{-1ex}%
\flushright\footnotesize\vspace{-4ex}%
$\mtho\longrightarrow$\\%
\vspace{-1ex}{#1}$\mtho$\vspace*{2ex}]%
{%\vspace{-1ex}%
\flushleft\footnotesize\vspace{-4ex}%
$\mtho\longleftarrow$\\%
\vspace{-1ex}{#1}$\mtho$\vspace*{2ex}}
}%

\nek{\pws} [1] {\cP(#1)}
\nek{\lla} {\,\land\,}

\nek{\snos} [1] {\,\footnote{\ #1}}
\nek{\snom}   {\,\footnotemark}
\nek{\snot} [1] {\footnotetext{\ #1}}

\nek{\ap} {{\hspace*{0.2ex}\boldsymbol\cdot\hspace*{0.2ex}}} 

\nek{\onto} {\overset{{\text{onto}}}{\longrightarrow}}

\nek{\La} {\Lambda}
\nek{\alo} {{\aleph_0}}

\nek{\leqv} {\,\eqv\,}

\nek{\dop} [1] {{#1}{}^{\complement}}

\nek{\ba}{\beta}

\nek{\etc} {{\sl etc}.}

\nek{\qand} {\quad\text{and}\quad}

\nek{\vt}{\vartheta}

\nek{\nse} {\om^{<\om}}
\nek{\bse} {2^{<\om}}

\nek{\namx} [1] 
{\overset{\text{\mtho$\hspace*{0.2ex}_\text{\Large\bf.}$}}{#1}}
\nek{\namy} [1] 
{\overset{\text{\mtho$\hspace*{0.5ex}_\text{\Large\bf.}$}}{#1}}
\nek{\namz} [1] 
{\overset{\text{\mtho$\hspace*{0.9ex}_\text{\Large\bf.}$}}{#1}}

\nek{\cD} {\mathscr D}

\nek{\np}{\newpage}

\nek{\sqf} {\sq^{\text{\tt fin}}}

\nek{\ka}{\kappa}

\nek{\pwor} [1] {{\text{\rm PWO}}(#1)}
\nek{\sepa} [1] {{\text{\ubf Sep}}(#1)}
\nek{\redu} [1] {{\text{\ubf Red}}(#1)}

\nek{\pet} {\text{\ubf PT}}
\nek{\ret} [1] {\res_{\hspace*{0.05ex}#1}}

\renek{\leq} {\leqslant}
\renek{\geq} {\geqslant}

\nek{\limp} {\mathrel{\,\imp\,}}

\nek{\spe} [1] {\text{\bf SC}(#1)}
\nek{\spf} [1] {\text{\bf SC}(#1)}
\nek{\spg} [1] {\text{\bf SC}^{<\om}(#1)}

\nek{\vys} [1] {\text{\tt hgt}(#1)} 

\nek{\cle} {\preccurlyeq}
\nek{\cls} {\prec}

\renek{\refname} {{\large\bf References}}

\nek{\dphi} {\mathbb\Phi}

\nek{\dU}{\mathbb U}
\nek{\uu}{^{\text{\tt fu}}}

\nek{\Tv}{\overrightarrow{T}}
\nek{\Sv}{\overrightarrow{S}}

\nek{\jta} {\boldsymbol\tau}
\nek{\jsg} {\boldsymbol\sg}
\nek{\jrho} {\boldsymbol\rho}
\nek{\jro} {\jrho}
\nek{\ju} {\boldsymbol u}
\nek{\jv} {\boldsymbol v}

\nek{\dn} {2^\om}
\nek{\bn} {\om^\om}

\nek{\dox}{{\namy {\boldsymbol x}}}

\nek{\rc} {\mathbf c}
\nek{\rd} {\mathbf d}
\nek{\rpi} {\dox}

\nek{\dplom} {\dd{\dP\lom}}
\nek{\dplo} {\dd\plo}

\nek{\plo} {\dP\lom}
\nek{\uplo} {(\dP\cup\dU)\lom}

\nek{\cll} {\mathrel{{\cls}\hspace*{-0.9ex}{\cls}}}

\nek{\vpj} [2] {\vpi^{#1}_{#2}}
\nek{\Phj} [2] {\Phi^{#1}_{#2}}

\nek{\vpt} [3] {\vpi^{#1}_{#2#3}}

\nek{\od} {\text{OD}}

\nek{\Eo} {\mathrel{{\text{\sf E}}_0}}

\nek{\zf} {\text{\ubf ZF}}

\nek{\uf} [2] {{\boldsymbol U}^\dphi_{#1}(#2)}
\nek{\ufi} [1] {{\boldsymbol U}^\dphi_{#1}}

\nek{\tf} [3] {{\boldsymbol U}^\dphi_{#1#2}(#3)}
\nek{\tfi} [2] {{\boldsymbol U}^\dphi_{#1#2}}

\nek{\tx} [2] {{\boldsymbol T}^\dphi_{#1}(#2)}
\nek{\txi} [1] {{\boldsymbol T}^\dphi_{#1}}

\nek{\ty} [3] {{\boldsymbol T}^\dphi_{#1#2}(#3)}

\nek{\lel} {\leq_{\rL}}

\renek{\rL} {\mathbf L}

\nek{\kc} [2] {C_{#1#2}}
\nek{\xc} [2] {C'_{#1#2}}
\nek{\kd} [2] {D_{#1#2}}
\nek{\kcp} [3] {C_{#1#2}^{#3}}
\nek{\kr} [3] {\jrho_{#1#2}^{#3}}
\nek{\kR} [2] {R_{#1#2}}

\nek{\jcp} [4] {C_{#1#2}^{#3#4}}
\nek{\jr} [4] {\jrho_{#1#2}^{#3#4}}

\renek{\imar}[1]{} \renek{\imae}[1] {}

\nek{\roo} [1] {\text{\tt stem}(#1)}

\nek{\cI} {\mathcal I}
\nek{\fss} [1] {\text{\bf FSS}(#1)}
\nek{\fsd} [1] {\text{\bf FSS}\lom(#1)}

\nek{\ptf} {\mathbb {PTF}}

\nek{\jp} {\mathbb p}
\nek{\bbu} {\mathbb u}
\nek{\bbp} {\mathbb p}
\nek{\bbq} {\mathbb q}
\nek{\bbpu} {\bbp{\lor}\bbu}

\nek{\mus} {multisystem}
\nek{\ms} [1] {\text{\ubf MS}(#1)}

\nek{\mut} {multitree}
\nek{\mt} [1] {\text{\ubf MT}(#1)}

\nek{\omim} {\omi^\cM}
\nek{\et}{\eta}             

\nek{\bas} {\underline s}
\nek{\bat} {\underline t}
\nek{\bah} {H}
\nek{\baH} {2^H}

\nek{\gh} {2^{H}} %{\mathfrak H}

\nek{\jvt} {\bbtheta}

\nek{\ivt} {\underline\vt^{-1}}

\nek{\wo} {\text{\ubf WO}}
\nek{\mto} {\mapsto}

\nek{\ahh} [1] {A_{#1}^\ast}
\nek{\bhh} [1] {B_{#1}^\ast}

\nek{\uH}{\mathbf h}

\nek{\bp} {\bar p}
\nek{\bq} {\bar q}

\nek{\mtp} [1] {\mt{#1}^+}

\nek{\Ord} {\text{\ubf Ord}}

\nek{\xh} {p}

\begin{document}

\title
{Counterexamples to countable-section $\ip12$ uniformization 
and $\varPi^1_3$ separation}

\author 
%, Lyubetsky]
{
%Ali~Enayat\thanks{address, email}  
%\and
Vladimir~Kanovei\thanks{IITP RAS and MIIT,
  Moscow, Russia, \ {\tt kanovei@googlemail.com} --- 
contact author. 
Partial support of 
RFFI grant 13-01-00006 acknowledged.}  
\and
Vassily~Lyubetsky\thanks{IITP RAS,
  Moscow, Russia, \ {\tt lyubetsk@iitp.ru} 
}}

\date 
{\today}

\maketitle

\begin{abstract}
We make use of a finite support product of 
the Jensen minimal $\varPi^1_2$ singleton forcing to 
define a model in which $\ip12$ Uniformization fails for a 
set with countable cross-sections. 
We also define appropriate submodels of the same model in 
which Separation fails for $\ip13$.
\end{abstract}

\parf{Introduction}
\las{cha1}

The uniformization problem, introduced by Luzin \cite{lhad}, 
is well known in modern set theory. 
(See Moschovakis~\cite{mDST} and Kechris~\cite{dst} for 
both older and more recent studies.)
In particular, it is known that every $\fs12$ set can be 
uniformized by a set of the same class $\fs12$, but on 
the other hand, there is a $\fp12$ set (in fact, a lightface 
$\ip12$ set), not uniformizable by any set in $\fp12$.

The negative part of this result cannot be strengthened much 
further in the direction of more complicated uniformizing 
sets since any $\fp12$ set admits a \dd{\fd13}uniformization 
assuming $\rV=\rL$ and admits a \dd{\fp13}uniformization 
assuming the existence of sharps (the Martin -- Solovay -- 
Mansfield theorem, \cite[8H.10]{mDST}).

However, the mentioned 
\dd{\fp12}non-uniformization theorem can be strengthened 
in the context of consistency. 
For instance, the $\ip12$ set 
$$
P=\ens{\ang{x,y}}{x,y\in\dn\land y\nin\rL[x]}
$$
is not uniformizable by any ROD (real-ordinal definable) set 
in the Solovay model and many other models of $\ZFC$ in which 
it is not true that $\rV=\rL[x]$ for a real $x$, and then 
the cross-sections of $P$ can be considered as ``large'', 
in particular, they are definitely uncountable.
Therefore one may ask:

\bvo
\lam{vo1}
Can such a ROD-non-uniformizable $\ip12$ set $P$ 
have the property 
that all his cross-sections are at most countable? 
\evo

This question is obviously connected with another 
question, initiated and briefly discussed at the 
{\it Mathoverflow\/} exchange desk\snos
{\label{snos1}
A question about ordinal definable real numbers. 
Mathoverflow, March 09, 2010. 
{\tt http://mathoverflow.net/questions/17608}. 
}
and at FOM\snos
{\label{snos2}%
Ali Enayat. Ordinal definable numbers. FOM Jul 23, 2010.
{\tt http://cs.nyu.edu/pipermail/fom/2010-July/014944.html}}
:

\bvo
\lam{vo2}
Is it consistent with $\ZFC$ that there is a 
\rit{countable} 
definable set of reals $X\ne\pu$ which has no OD 
(ordinal definable) elements.
\evo

Ali Enayat (Footnote~\ref{snos2}) conjectured that 
Question~\ref{vo2} can be solved in the positive by the   
finite-support product $\plo$ 
of countably many copies of the Jensen ``minimal $\ip12$ real 
singleton forcing'' $\dP$ defined in \cite{jenmin} 
(see also Section 28A of \cite{jechmill}).
Enayat demonstrated that a symmetric part of the \dd\plo generic 
extension of $\rL$ definitely yields a model of $\zf$ 
(not a model of $\ZFC$!) 
in which there is a Dedekind-finite infinite \od\ set of 
reals with no \od\ elements. 

Following the mentioned conjecture, we proved in \cite{klOD} 
that indeed it is true in a \dd\plo generic extension of $\rL$  
that the set of \dd\dP generic reals is a countable non-empty  
$\ip12$ set with no $\od$ elements.\snos  
{We also proved in \cite{klE} that the existence of a 
$\ip12$ \dd\Eo class with no $\od$ elements is consistent 
with $\ZFC$, using a \dd\Eo invariant version of the Jensen 
forcing.}  
Using a finite-support product $\prod_{\xi<\omi}{\dP_\xi}\lom$, 
where all $\dP_\xi$ are forcings similar to, 
but different from, Jensen's forcing $\dP$ (and from each other), 
we answer Question~\ref{vo1} in the positive.

\bte
\lam{Tun}
In a suitable generic extension of\/ $\rL$, 
it is true that there is a lightface\/ 
$\ip12$ set\/ $P\sq\dn\ti\dn$ whose all cross-sections\/ 
$P_x=\ens{y}{\ang{x,y}\in P}$ are at most countable, but\/ 
$P$ is not uniformizable by a ROD set.
\ete

Using an appropriate generic extension of a 
submodel of the same model, similar to  
models considered in Harrington's unpublished notes 
\cite{h74}, we also prove 

\bte
\lam{Tsep}
In a suitable generic extension of\/ $\rL$, 
it is true  that there is a pair of 
disjoint lightface\/ $\ip13$ sets\/ $X,Y\sq\dn$, not 
separable by disjoint\/ $\fs13$ sets, and hence\/ 
$\fp13$ Separation and\/ 
$\ip13$ Separation fail.
\ete

This result was first proved by Harrington in \cite{h74} 
on the base of almost disjoint forcing of Jensen -- 
Solovay \cite{jsad}, 
and in this form 
has never been published, but was mentioned, \eg, in 
\cite[5B.3]{mDST} and \cite[page 230]{hin}.
A complicated alternative proof of Theorem~\ref{Tsep} 
can be obtained 
%in \cite{k83} 
with the help of \rit{countable-support} products and 
iterations of Jensen's forcing studied earlier in 
\cite{abr,k79,k83}. 
The \rit{finite-support} approach which we pursue here yields 
a significantly more compact proof. 
As far as Theorem~\ref{Tun} is concerned, countable-support 
products and iterations hardly can lead to the 
countable-section non-uniformization results.

We recall that the $\fp13$ Separation 
\rit{hold} in $\rL$, the constructible universe. 
Thus Theorem~\ref{Tsep} in fact shows that the $\fp13$ 
Separation principle is ``killed'' in an appropriate generic 
extension of $\rL$.  
It would be interesting to find a generic extension in which, 
the other way around, the $\fs13$ Separation 
(false in $\rL$) holds.

%See an application of \dd{\om_2}long 
%iterations of Jensen's forcing in \cite{abr}. 

\back
The authors thank Jindra Zapletal and Ali Enayat for 
fruitful discussions.
\eack

%\cite{abr2}

%\newpage

\vyk{
\snos
{\label{nou}%
The most elementary counterexample can be defined as follows. 
We begin with a pair of disjoint $\is12$ sets $X,Y\sq\bn,$ 
nonseparable by a $\fd12$ set. 
Let $P$ consist of all pairs $\ang{x,k}\in\bn\ti\om$ such 
that $x\nin X\land k=0$ or $x\nin Y\land k=1$; this $P$ 
is $\ip12$. 
If $Q\sq P$ is a $\fp12$ set which uniformizes $P$ then 
$X'=\ens{x}{\ang{x,1}\in Q}$ and $Y'=\ens{x}{\ang{x,0}\in Q}$ 
are mutually complimentary $\fp12$ sets, hence $X'$ is a 
$\fd12$ set which separates $X$ from $Y$, a contradiction.}
}

\parf{Trees, perfect-tree forcing notions, splitting} 
\las{tre}
\label{ptf}

Let $\bse$ be the set of all strings (finite sequences) 
of numbers $0,1$.
If $t\in\bse$ and $i=0,1$ then 
$t\we k$ is the extension of $t$ by $k$. 
If $s,t\in\bse$ then $s\sq t$ means that $t$ extends $s$, while 
$s\su t$ means proper extension. 
If $s\in\bse$ then $\lh s$ is the length of $s$,  
and $2^n=\ens{s\in\bse}{\lh s=n}$ (strings of length $n$).%

A set $T\sq\bse$ is a \rit{tree} iff 
\index{tree}%
%it is an initial segment, that is, 
for any strings $s\su t$ in $\bse$, if $t\in T$ then $s\in T$.
Thus every non-empty tree $T\sq\bse$ contains the empty 
string $\La$. 
If $T\sq\bse$ is a tree and $s\in T$ then put 
$T\ret s=\ens{t\in T}{s\sq t\lor t\sq s}$. 
%this is a tree as well.

Let $\pet$ be the set of all \rit{perfect} trees 
$\pu\ne T\sq \bse$. 
\imar{pet}%
Thus a non-empty tree $T\sq\bse$ belongs to $\pet$ iff 
it has no endpoints and no isolated branches. 
Then there is a largest string $s\in T$ such that 
$T=T\ret s$; it is denoted by $s=\roo T$   
(the {\it stem\/} of $T$); 
we have $s\we 1\in T$ and $s\we 0\in T$ in this case.

Each perfect tree $T\in\pet$ defines  
$ 
[T]=\ens{a\in\dn}{\kaz n\,(a\res n\in T)}\sq\dn
$, 
the perfect set of all \rit{paths through $T$}. 

\bdf
\lam{ptf2}
A {\ubf perfect-tree forcing notion} is any set 
$\dP\sq\pet$ such that if $u\in T\in\dP$ then $T\ret u\in \dP$.
Let $\ptf$ be the set of all such  
$\dP\sq\pet$.
\edf

Such a set $\dP$ can be considered as a forcing notion 
(if $T\sq T'$ then $T$ is a stronger condition); 
such a forcing $\dP$ adds a real in $\dn$. 

\bpri
\lam{p0}
If $s\in\bse$ then the 
%``Baire interval'' 
tree $I_s=\ens{t\in\bse}{s\sq t\lor t\sq s}$ 
belongs to $\pet$ and the set $\dP_0=\ens{I_s}{s\in\bse}$ 
is a perfect-tree forcing.
\epri

\ble
\lam{suz}
If\/ $\dP,\dP'\in\pet$, $T\in\dP$, $T'\in\dP'$, then 
there are trees\/ $S\in\dP$, $S'\in\dP'$ such that\/ 
$S\sq T$, $S'\sq T'$, and\/ $[S]\cap[S']=\pu$.
\ele
\bpf
If $T=T'$ then let $s=\roo T$ and $S=T\ret{s\we0}$, 
$S'=T'\ret{s\we1}$.
If say $T\not\sq T'$ then let $s\in T\bez T'$, 
$S=T\ret{s}$, and simply $S'=T'$.
\epf

%\parf{Splitting construction over a perfect set forcing} 
%\las{spe}

%\bdf
%\lam{sped}
If $\dP\in\ptf$ then let $\fss \dP$ be the set of all 
\rit{finite splitting systems} over $\dP$, that is,  
systems of the form 
$\vpi=\sis{T_s}{s\in2^{< n}}$, where 
$n=\vys\vpi<\om$ (the height of $\vpi$), 
each value $T_s=\vpi(s)$ is a tree in $\dP$, and
\ben
\fenu
\itla{spe2}
if $s\we i\in 2^{<n}$ ($i=0,1$) then $T_{s\we i}\sq T_s$ and  
$\roo{T_{s}}\we i\sq\roo{T_{s\we i}}$ --- 
\imar{spe2}
it easily follows that 
$[T_{s\we0}]\cap [T_{s\we0}]=\pu$.
%\qed
\een
%\eDf 
% 
Let $\vpi,\psi$ be systems in $\fss\dP$.
Say that 
\bit
\item[$-$]
$\vpi$ \rit{extends} $\psi$, symbolically $\psi\cle\vpi$, if 
$n=\vys\psi\le\vys\vpi$ and $\psi(s)=\vpi(s)$ for 
all $s\in2^{<n}$;

\item[$-$]
\rit{properly extends} $\psi$, 
symbolically $\psi\cls\vpi$, if in 
addition $\vys\psi<\vys\vpi$;

\item[$-$]
\rit{reduces} $\psi$, if $n=\vys\psi=\vys\vpi$, 
$\vpi(s)\sq\psi(s)$ for all $s\in 2^{\vys\vpi-1},$ and 
$\vpi(s)=\psi(s)$ for all $s\in 2^{<\vys\vpi-1}$.
\eit
In other words, reduction allows to shrink trees in the top 
layer of the system, but keeps intact those in the lower 
layers.

The empty system $\La$ is the only one in 
$\fss\dP$ satisfying $\vys\La=0$.
To get a system $\vpi\in\fss\dP$ with $\vys\vpi=1$ take any 
$T\in\dP$ and put $\vpi(\La)=T$.
The next lemma provides systems of bigger height.

\ble
\lam{n+1}
Assume that\/ $\dP\in\ptf$.
If\/ $n\ge1$ and\/ $\psi=\sis{T_s}{s\in2^{< n}}\in\fss\dP$ 
then there is a system\/ 
$\vpi=\sis{T_s}{s\in2^{< n+1}}\in\fss\dP$ which properly extends
\/ $\psi$.
\ele
\bpf
If $s\in2^{n-1}$ and $i=0,1$ then let 
$T_{s\we i}=T_s\ret{\roo{T_s}\we i}$. 
\epf

\bcor
\lam{infty}
Let\/ $\dP\in\ptf$.
Then there is an\/ \dd\cls in\-creas\-ing sequence\/ 
$\sis{\vpi_n}{n<\om}$ of systems in\/ $\fss\dP$. 
In this case the limit system\/ 
$\vpi=\bigcup_n\vpi_n=\sis{T_s}{s\in\bse}$
satisfies\/ \ref{spe2} of Section~\ref{tre}
%Definition~\ref{sped} 
on the whole domain\/ $\bse,$ 
$T=\bigcap_n\bigcup_{s\in2^n}T_s$ is 
a perfect tree in\/ $\pet$ 
{\rm(yet not necessarily in $\dP$)}, 
and\/ $[T]=\bigcap_n\bigcup_{s\in2^n}[T_s]$. \qed
\ecor

Say that a tree $T$ \rit{occurs in\/ $\vpi\in\fss\dP$} if 
$T=\vpi(s)$ for some $s\in 2^{<\vys\vpi}$.

\parf{Multitrees and splitting multisystems} 
\las{mul}

Suppose that $\vt\in\Ord$ and $\bbp=\sis{\dP_\xi}{\xi<\vt}$ 
is a sequence of sets $\dP_\xi\in\ptf$. 
We'll systematically consider such sequences below, and 
if $\bbq=\sis{\dQ_\xi}{\xi<\vt}$ is another such a sequence 
of the same length then let 
$\bbp{\lor}\bbq  =\sis{\dP_\xi\cup\dQ_\xi}{\xi<\vt}$.

\bdf
\lam{muts}
A \rit{\dd\bbp\mut} is a ``matrix'' of the form 
$\jta=\sid{T_{\xi k}}{\xi<\vt}{k<\om}$, where each 
$\jta(\xi,k)=T_{\xi k}$ belongs to $\dP_\xi$, and the 
\rit{support} $\abs\jta=\ens{\ang{\xi,k}}{T_{\xi k}\ne\bse}$ 
is finite.
Let $\mt\bbp$ be the set  of all \dd\bbp\mut s.
If $\jta\in\mt\bbp$ then let
$$
[\jta]=\ens{x\in2^{\vt\ti\om}}
{\kaz \ang{\xi,k}\in\abs\jta\,(x(\xi,k)\in [\jta(\xi,k)])}\,;
$$
this is a cofinite-dimensional perfect cube in $2^{\vt\ti\om}$.

A \rit{\dd\bbp\mus} is a ``matrix'' of the form 
$\Phi=\sid{\vpi_{\xi m}}{\xi<\vt}{ m<\om}$, where each 
$\Phi(\xi, m)=\vpi_{\xi m}$ belongs to $\fss{\dP_\xi}$, and the 
\rit{support} $\abs\Phi=\ens{\ang{\xi, m}}{\vpi_{\xi m}\ne\bse}$ 
is finite.
Let $\ms\bbp$ be the set  of all \dd\bbp\mus s.

Say that a \mut\ $\jta=\sid{T_{\xi k}}{\xi<\vt}{k<\om}$ 
\rit{occurs} in a \mus\ 
$\Phi=\sid{\vpi_{\xi m}}{\xi<\vt}{ m<\om}$  
if $\abs\jta\sq\abs\Phi$ and for each $\ang{\xi,k}\in\abs\jta$ 
there is a number $ m<\om$ and a 
string $s\in\bse$ with $\lh s<\vys{\vpi_{\xi m}}$ 
such that $T_{\xi k}=\vpi_{\xi m}(s)$.
\edf
 
The set $\mt\bbp$ is equal to the finite support product 
$\prod_{\xi<\vt}(\dP_\xi)^\om$ 
of $\vt\ti\om$-many factors, with each factor $\dP_\xi$ 
in \dd\om many copies. 
Accordingly, the set $\ms\bbp$ is equal to the 
finite support product 
$\prod_{\xi<\vt}(\fss{\dP_\xi})^\om$ 
of $(\vt\ti\om)$-many factors, with each factor $\fss{\dP_\xi}$ 
in \dd\om many copies. 
We order $\mt\bbp$ componentwise: $\jsg\leq\jta$ 
%($\jsg$ is stronger) 
iff $\jsg(\xi,k)\sq\jta(\xi,k)$ in $\dP_\xi$ for all $\xi,k$.
The forcing $\mt\bbp$ adds a ``matrix''   
$\sid{x_{\xi k}}{\xi<\vt}{k<\om}$, where each $x_{\xi k}\in\dn$ 
is a \dd{\dP_\xi}generic real. 

If $\Phi,\Psi\in\ms\bbp$ then we define 
\bit
\item[$-$] 
$\Psi\cle\Phi$ \ iff \ $\Psi(\xi, m)\cle\Phi(\xi, m)$ 
(in $\fss{\dP_\xi}$) for all $\xi, m$;

\item[$-$] 
$\Phi$ reduces $\Psi$ \ iff \ $\abs\Psi\sq\abs\Phi$ and 
$\Phi(\xi,m)$ reduces $\Psi(\xi,m)$ for all pairs 
$\ang{\xi, m}\in\abs\Psi$;

\item[$-$] 
$\Phi\cll\Psi$ \ iff \  $\abs\Phi\sq\abs\Psi$ 
and $\Phi(\xi, m)\cls\Psi(\xi, m)$ for all 
$\ang{\xi, m}\in\abs\Phi$.
\eit 
%Then $\Psi\cls\Phi$ means that $\Psi\cle\Phi$ and 
%$\Psi(\xi, m)\cls\Phi(\xi, m)$ for at least 
%one pair $\ang{\xi, m}$. 

\ble
\lam{xr}
If $\Phi\cll\Psi$ and $\Phi'$ reduces $\Psi$ then still 
$\Phi\cll\Phi'$ and $\Phi\cle\Phi'$.\qed
\ele

\parf{Jensen's extension of a perfect tree forcing} 
\las{jex}

Let $\zfc'$ be the subtheory of $\zfc$ including all 
axioms except for the 
power set axiom, plus the axiom saying that $\pws\om$ exists. 
(Then $\omi$ and continual sets like $\pet$ exist as well.)
Let $\cM$ be a countable transitive model of $\ZFC'$. 

Suppose that $\bbp=\sis{\dP_\xi}{\xi<\jvt}\in\cM$ 
is a sequence of (countable) sets $\dP_\xi\in\ptf$, 
of length $\jvt<\omim$.  
Then the sets ${\dP_\xi}$ and $\fss{\dP_\xi}$ for all $\xi<\jvt$, 
as well as the sets $\mt\bbp$ and $\ms\bbp$, belong to $\cM$, too. 

\bdf
\lam{dPhi}
(i) 
Let us fix any \dd\cle increasing sequence 
$\dphi=\sis{\Phi^j}{j<\om}$ of \mus s 
$\Phi^j=\sid{\vpt j\xi m}{\xi<\jvt}{ m<\om}\in\ms\bbp$, 
\rit{generic over\/ $\cM$} in the sense that it intersects 
every set $D\in\cM\yd D\sq\ms\bbp$, dense in $\ms\bbp$\snos
{Meaning that for any $\Psi\in\ms\bbp$ there is $\Phi\in D$ 
with $\Psi\cle\Phi$.}% 
.

(ii) 
Suppose that $\xi<\jvt$ and $ m<\om$. 
In particular, the sequence $\dphi$ intersects every set 
of the form   
$$
D_{\xi m h}=\ens{\Phi\in\ms\bbp}
{\vys{\Phi(\xi, m)}\ge h}\,,\quad\text{where }h<\om\,. 
$$
It follows that the sequence 
$\sis{\vpt j\xi m}{j<\om}$ of systems 
$\vpt j\xi m\in\fss{\dP_\xi}$ 
%is \rit{eventually strictly increasing}, so that 
satisfies 
$\vpt j\xi m\cls \vpt {j+1}\xi m$ 
for infinitely many indices $j$  
(and $\vpt j\xi m=\vpt {j+1}\xi m$ for other $j$).

(iii)
We conclude that the limit system 
$\vpt\iy\xi m=\bigcup_{j<\om}\vpt j\xi m$ has the form 
$\sis{\ty \xi m s}{s\in\bse}$ such that each $\ty \xi m s$
is a tree in $\dP_\xi$, and if $j<\om$ then we have 
$\vpt j\xi m=\sis{\ty \xi m s}{s\in 2^{<h(j,\xi, m)}}$, 
where $h(j,\xi, m)=\vys{\vpt j\xi m}$. 

(iv)
Moreover, by Corollary \ref{infty}, the trees 
$$
\textstyle
\tfi\xi m=\bigcap_n\bigcup_{s\in2^n}\ty\xi m s\,,\quad  
\tf\xi m s=\bigcap_{n\ge \lh s}
\bigcup_{t\in2^n,\:s\sq t}\ty\xi m t
$$   
belong to $\pet$ (not necessarily to $\dP_\xi$) 
for each $s\in\bse;$ 
thus $\tfi\xi m=\tf\xi m \La$. 
%In fact $\tf\xi m s=\tfi\xi m\cap \ty\xi m s$ by 
%\ref{spe2} of Section~\ref{tre}. 
%Definition \ref{sped}.

(v) 
If $\xi<\jvt$ then let 
$\dU_\xi=\ens{\tf\xi m s}{ m<\om\land s\in\bse}$. 

Let $\bbu=\sis{\dU_\xi}{\xi<\jvt}$. 

Finally let  
$\bbpu=\sis{\dP_\xi\cup\dU_\xi}{\xi<\jvt}$.
\edf

\ble
\label{disj}
\ben
\renu
\itla{disj1}
if\/ $\ang{\xi,m}\ne\ang{\et,n}$ then\/ 
$[\tfi\xi m]\cap[\tfi\et n]=\pu\;;$

\itla{disj2}
if\/ $\xi<\jvt$, $m<\om$, $s\in\bse,$ then\/ 
$\tf\xi m s=\tfi\xi m\cap \ty\xi m s\;;$

\itla{disj3}
if\/ $\xi<\jvt$, $m<\om$, and 
strings\/ $s\sq t$ belong to\/ $\bse$ 
then\/ 
$[\ty\xi ms]\sq[\ty\xi mt]$ and\/ $[\tf\xi ms]\sq[\tf\xi mt]\;;$

\itla{disj4}
If\/ $\xi<\jvt$, $m<\om$, and strings\/ $t'\ne t$ in\/ $\bse$ 
are\/ \dd\sq incomparable then\/ 
$[\tf\xi m{t'}]\cap[\tf\xi mt]=[\ty\xi m{t'}]
\cap[\ty\xi mt]=\pu$. 
\een 
\ele
\bpf
\ref{disj1} 
By Lemma~\ref{suz}, the set $D$ of all \mus s $\Phi$ such that
the pairs $\ang{\xi,m}\yi\ang{\et,n}$ belong to $\abs\Phi$ and,
for some $h<\min\ans{\vys{\Phi(\xi,m)},\vys{\Phi(\et,n)}}$, 
we have $[\Phi(\xi,m)(s)]\cap [\Phi(\et,n)(t)]=\pu$  for all 
$s,t\in 2^h$, is dense.

\ref{disj2} 
easily follows from \ref{spe2} of Section~\ref{tre}.
\ref{disj3} 
is obvious.

\ref{disj4} 
Note that $[\vpi(s\we0)]\cap[\vpi(s\we1)]=\pu$ for any system 
$\vpi$ with $\vys\vpi>1+\lh s$ by \ref{spe2} of 
Section~\ref{tre}.
Therefore $[\ty\xi m{s\we 0}]\cap[\ty\xi m{s\we 1}]=\pu$.
\epf

It follows that if $U\in\bigcup_{\xi<\jvt}\dU_\xi$ then 
there is a unique triple of $\xi<\jvt$, $m<\om$, and 
$s\in\bse$ such that $U=\tf\xi m{s}$!

\ble
\lam{uu1}
If\/ $\xi<\jvt$ then  
the sets\/ $\dU_\xi$ and\/ 
$\dP_\xi\cup\dU_\xi$ belong to\/ $\ptf$.
% of Section~\ref{ptf}.
\qed
\ele

\ble
\lam{uu2}
Let\/ $\xi<\jvt$. 
The set\/ $\dU_\xi$ is dense in\/ $\dU_\xi\cup\dP_\xi$. 
\ele
\bpf
If $T\in\dP_\xi$ then 
the set $D(T)$ of all \mus s   
$\Phi=\sid{\vpi_{\xi m}}{\xi<\jvt}{ m<\om}$ in $\ms\bbp$,
such that $\vpi_{\xi m}(\La)=T$ for some $k$, belongs to $\cM$ 
and obviously is dense in $\ms\bbp$. 
It follows that $\Phi^J\in D(T)$ for some $J<\om$, 
by the choice of $\dphi$. 
Then $\ty\xi m\La=T$ for some $ m<\om$. 
However $\tf\xi m\La\sq \ty\xi m\La$. 
\epf

\ble
\lam{uu3}
If\/ $\xi<\jvt$ and a set\/ $D\in\cM$, $D\sq\dP_\xi$ 
is pre-dense in\/ $\dP_\xi$, 
and\/ $U\in\dU_\xi$, then\/ $U\sqf\bigcup D$, 
that is, there is a finite set\/ $D'\sq D$ with
$U\sq\bigcup D'$. 
\ele
\bpf
Suppose that $U=\tf\xi Ms$, $M<\om$ and $s\in\bse.$ 
Consider the set $\Da\in\cM$ of all \mus s   
$\Phi=\sid{\vpi_{\xi  m}}{}{}\in \ms\bbp$ such that 
$\ang{\xi,M}\in\abs\Phi$, $\lh s<h=\vys{\vpi_{\xi M}}$, 
and for each $t\in 2^{h-1}$ there is a tree $S_t\in D$ with  
$\vpi_{\xi M}(t)\sq S$.
The set $\Da$ 
is dense in $\spg\dP$ by the pre-density of $D$. 
Therefore there is an index $J$ such that $\Phi^J$ belongs 
to $\Da$.
Let this be witnessed by trees $S_t\in D\yt t\in 2^{h-1},$ 
where $\lh s<h=\vys{\vpt J\xi M}$, so that $\vpt J\xi M(t)\sq S_t$. 
Then 
$$
\textstyle
U=\tf\xi Ms\sq\tf\xi M\La\sq\bigcup_{t\in 2^{h-1}}\vpt J\xi M(t)
\sq\bigcup_{t\in 2^{h-1}}S_t\sq\bigcup D'
$$ 
by construction, where $D'=\ens{S_t}{t\in 2^{h-1}}\sq D$ 
is finite.
\epf

\ble
\lam{uu4}
If\/ 
%$\xi<\jvt$ and 
a set\/ $D\in\cM,$ $D\sq\mt\bbp$  
is pre-dense in\/ $\mt\bbp$ then it remains pre-dense in\/ 
$\mt\bbpu$. 
\ele
\bpf
Given a \mut\ $\jta\in\mt\bbpu$,  prove that 
$\jta$ is compatible in $\mt\bbpu$ with a \mut\ 
$\jsg\in D$.
For the sake of brevity, assume that $\jta\in\mt\bbu$ and 
$\abs\jta=\ans{\ang{\eta,K},\ang{\za,L}}$,  
where $\za<\eta<\jvt$ and $K,L<\om$. 
Then by construction 
$\jta(\eta,K)=\tf\eta{M}s$ 
%=\tfi\eta K\cap\ty\eta Ks$ 
and 
$\jta(\za,L)=\tf\za{N}t$ 
%=\tfi\za L\cap\ty\za Lt$  
for some $M,N<\om$ and $s,t\in\bse$.

Consider the set $\Da\in\cM$ of all \mus s   
$\Phi=\sid{\vpi_{\xi m}}{\xi<\jvt}{m<\om}\in \ms\bbp$ 
such that there are strings $s',t'\in\bse$ 
with $s\su s'$, $t\su t'$, 
$\lh s'<\vys{\vpi_{\eta M}}$, $\lh t'<\vys{\vpi_{\za N}}$, 
and \mut s\ $\jsg\in D$ and $\jsg'\in\mt\bbp$, such that 
$\jsg'\leq\jsg$ and $\jsg'$ occurs in $\Phi$ in such a way that 
$\jsg'(\eta,K)=\vpi_{\et M}(s')$ and 
$\jsg'(\za,L)=\vpi_{\za N}(t')$.
%
%$\vpi_k(s')\sq U\cap S$, $\vpi_\ell(t')\sq V\cap T$.

The set $\Da$ 
is dense in $\ms\bbp$ by the pre-density of $D$. 
Therefore there is an index $j$ such that $\Phi^j$ belongs 
to $\Da$.
Let this be witnessed by strings $s',t'\in\bse,$ and \mut s 
$\jsg\in D$, and $\jsg'\in\mt\bbp\yt \jsg'\leq\jsg$, as above. 
In other words, $s\su s'$, $t\su t'$, 
$\lh s'<\vys{\vpt j{\eta}{M}}$, 
$\lh t'<\vys{\vpt j{\za}{N}}$,
and $\jsg'$ occurs in $\Phi$ in such a way that 
$\jsg'(\eta,K)=\vpt j{\eta}{M}(s')$ and 
$\jsg'(\za,L)=\vpt j{\za}{N}(t')$.
The set 
$\abs{\jsg'}=\ans{\ang{\xi_1,k_1},\ang{\xi_2,k_2},\dots,
\ang{\xi_n,k_n}}\sq\jvt\ti\om$ 
is finite and contains the pairs $\ang{\eta,K},\ang{\za,L}$; 
let, say, $\ang{\xi_1,k_1}=\ang{\eta,K}$, 
$\ang{\xi_2,k_2}=\ang{\za,L}$.

And if $i=1,2,\dots,n$ then by definition  
$\jsg'(\xi_i,k_i)=\vpt j{\xi_i}{m_i}(s_i)=\ty {\xi_i}{m_i}{s_i}$ 
holds for some $m_i<\om$ and $s_i\in\bse.$
In particular 
$\jsg'(\eta,K)=\vpt j{\eta}{M}(s')=\ty {\eta}{M}{s'}$ and 
$\jsg'(\za,L)=\vpt j{\za}{N}(t')=\ty {\za}{N}{t'}$, 
for $i=1,2$. 

Consider the \mut\ $\jta'\in\mt\bbu$ defined so that 
$\abs{\jta'}=\abs{\jsg'}$ and 
$\jta'(\xi_i,k_i)=\tf{\xi_i}{m_i}{s_i}$
for all $i=1,\dots,n$. 
In particular 
$\jta'(\eta,K)=\tf{\eta}{M}{s'}$ and 
$\jta'(\za,L)=\tf{\za}{N}{t'}$.
Then $\jta'\leq\jsg'$ 
(since $\tf{\xi_i}{m_i}{s_i}\sq \ty{\xi_i}{m_i}{s_i}$), 
therefore $\jta'\leq\jsg\in D$. 

It remains to prove that $\jta'\leq\jta$, which amounts to 
$\jta'(\eta,K)\sq\jta(\eta,K)$ and $\jta'(\za,L)\sq\jta(\za,L)$.
However $\jta(\eta,K)=\tf\eta M s\sq \tf\eta M{s'}=\jta'(\eta,K)$
since $s\su s'$, and the same for the pair $\ang{\za,L}$.
\epf

\vyk{
\ble
\lam{le:uu}
In the assumptions above, the set\/ $(\dP\cup\dU)\uu$ 
of all finite unions of trees in $\dP\cup\dU$ is still a 
perfect-tree forcing.\qed 
\ele
}

\parf{Forcing a real away of a pre-dense set} 
\las{saway}

Let $\cM$ be still a countable transitive model of $\ZFC'$ 
and $\bbp=\sis{\dP_\xi}{\xi<\omi^\cM}\in\cM$ be as in 
Section~\ref{jex}. 
The goal of the following Theorem~\ref{K} is to prove that, 
under the conditions and notation of Definition~\ref{dPhi}, 
if $\xi<\jvt$ and $c$ is a \dd{\mt\bbp}name of a real in $\dn$ 
then 
it is forced by the extended forcing $\mt\bbpu$ 
that $c$ does not belong to sets $[U]$ where 
$U$ is a tree in $\dU_\xi$ --- unless $c$ is a name 
of one of generic reals $x_{\xi k}$ themselves.
We begin with a suitable notation.

\bdf
\lam{rk}
A \rit{\dd{\mt\bbp}real name} is a system 
$\rc=\sis{\kc ni}{n<\om,\, i<2}$ of sets $\kc ni\sq\mt\bbp$ 
such that each set $C_n=\kc n0\cup \kc n1$ is 
dense or at least pre-dense in $\mt\bbp$ 
and if $\jsg\in \kc n0$ and $\jta\in \kc n1$ then $\jsg,\jta$ are 
incompatible in $\mt\bbp$.

If a set $G\sq\mt\bbp$ is \dd{\mt\bbp}generic at least over 
the collection of all  sets $C_n$ then we define 
$\rc[G]\in\dn$ so that $\rc[G](n)=i$ iff $G\cap \kc ni\ne\pu$.
\edf

Thus any \dd{\mt\bbp}real name $\rc=\sis{\kc ni}{}$ 
is a \dd{\mt\bbp}name for a real in $\dn.$ 

Recall that ${\mt\bbp}$ adds a generic sequence 
$\sis{x_{\xi k}}{\xi<\jvt,k<\om}$ of reals $x_{\xi k}\in\dn.$ 

\bpri
\lam{proj}
If $\xi<\jvt$ and $k<\om$ then  
define a \dd{\mt\bbp}real name 
$\rpi_{\xi k}=\sis{\jcp ni\xi k}{n<\om,\,i<2}$
such that each set $\jcp ni\xi k$ contains a single \mut\ 
$\jr ni\xi k\in\mt\bbp$, such that 
$\abs{\jr ni\xi k}=\ans{\ang{\xi,k}}$ and finally   
$\jr ni\xi k(\xi,k)=\kR ni$, where 
$$
\kR ni=\ens{s\in\bse}{\lh s>n\imp s(n)=i}\,.
$$
Then $\rpi_{\xi k}$ is a \dd{\mt\bbp}real name 
of the real $x_{\xi k}$, the $(\xi,k)$th 
term of a \dd{\mt\bbp}generic sequence 
$\sis{x_{\xi k}}{\xi<\jvt,\,k<\om}$. 
\epri

Let $\rc=\sis{\kc ni}{}$ and $\rd=\sis{\kd ni}{}$ 
be \dd{\mt\bbp}real names. 
Say that  $\jta\in \mt\bbp$:
\bit
\item
\rit{directly forces\/ $\rc(n)=i$}, 
where $n<\om$ and $i=0,1$, iff there is a finite set 
$\Sg\sq \kc ni$ such that $[\jta]\sq\bigcup_{\jsg\in\Sg}[\jsg]$; 
%(that is, the tree $T=\jta(k)\in\pet$ satisfies 
%$x(n)=i$ for all $x\in[T]$); 

\item
\rit{directly forces\/ $s\su\rc$},  
where $s\in\bse,$ iff for all $n<\lh s$, $\jta$ 
directly forces $\rc(n)=i$, where $i=s(n)$; 

\item
\rit{directly forces\/ $\rd\ne\rc$}, iff there are strings 
$s,t\in\bse,$ incomparable in $\bse$ and such that  
$\jta$ directly forces $s\su\rc$ 
and $t\su\rd$; 

\item
\rit{directly forces\/ $\rc\nin[T]$},  
where $T\in\pet$, iff there is a string $s\in\bse\bez T$ 
such that $\jta$ directly forces $s\su \rc$; 
\vyk{
\item
\,{\ubf[applicable only for $\jta\in\plo$]}
\rit{weakly\/ \dd\plo forces\/ $\rc\nin[T]$},  
iff the set of all \mut s 
$\jsg\in\plo$ that directly force\/ $\rc\nin[T]$ 
is dense in $\plo$ below $\jta$.
} 
\eit

\bte
\lam{K}
In the assumptions of Definition~\ref{dPhi}, suppose that\/ 
$\et<\vt$, $\rc= 
\sis{C_m^i}{m<\om,\,i<2}\in\cM$ is a\/ \dd{\mt\bbp}real name, 
and for all\/ $k$ the set
$$
D(k)=\ens{\jta\in\mt\bbp}{\jta\text{ directly forces }
\rc\ne\rpi_{\et k}}
$$
is dense in\/ $\mt\bbp$. 
Let\/ $\ju\in\mt\bbpu$, $\eta<\jvt$, and\/ $U\in \dU_\et$.
Then there is\/ 
a stronger \mut\/ $\jv\in\mt\bbu\yd \jv\leq\ju$, which 
directly forces\/ $\rc\nin[U]$.
\ete

%\vyk{
\bpf
By construction $U\sq \tfi\et M$ for some $M<\om$; 
thus we can assume that simply $U=\tfi\et M$. 
The indices 
$\et\yi M$ are fixed in the proof. 
We can assume by Lemma~\ref{uu2} that $\ju\in\mt\bbu$. 
The support 
$\abs\ju=\ans{\ang{\xi_1,k_1},\dots,\ang{\xi_\nu,k_\nu}}
\sq\jvt\ti\om$
is a finite set ($\nu<\om$), 
and if $i=1,\dots,\nu$ then, as $\ju\in\mt\bbu$, 
there is a string $s_i$ and a number $m_i$ such that 
$\ju(\xi_i,k_i)=\tf{\xi_i}{m_i}{s_i}$.
We can assume that 
\ben
\aenu
\itla{ass2}
if $i\ne i'$ and $\xi_i=\xi_{i'}$ then $k_i\ne k_i'$;
\imar{ass2}

\itla{ass3}
$s_i\ne s_{i'}$ whenever $i\ne i'$, 
\imar{ass3}
and there is $h<\om$ 
such that $\lh{s_i}=h\yd\kaz i\;;$\snos
{If $s_i\su s'_i\in\bse$ for all $i$, and $\ju'\in\mt\bbu$, 
$\abs{\ju'}=\abs{\ju}$ and 
$\ju'(\xi_i,k_i)=\tf{\xi_i}{m_i}{s'_i}$ for all $i$, 
then $\ju'\leq \ju$.
Thus if we prove the theorem for $\ju'$ then it implies 
the result for $\ju$ as well. 
}

\vyk{
\itla{ass1}
let $H=2^h$ --- then 
$\nu\ge H$, 
\imar{ass1}
$\xi_1=\dots=\xi_{H}=\et$, $m_1=\dots=m_{H}=M$, and 
$\ans{s_1,\dots,s_H}=\ens{s\in\bse}{\lh s=h}$;

\itla{ass4}
if $H<i\le\nu$ then 
\imar{ass4}
$\ang{\xi_i,m_i}\ne \ang{\et,M}$.
}%
\itla{ass1}
there is a number $\mu\le\nu$ 
such that $\xi_1=\dots=\xi_\mu=\et$ and 
\imar{ass1}
$m_1=\dots=m_\mu=M$ (then $\mu\le 2^h$), 
but if $\mu<i\le\nu$ then 
$\ang{\xi_i,m_i}\ne \ang{\et,M}$.
\een
In these assumptions, define a \mut\ $\jta\in\mt\bbp$ so 
that $\abs\jta=\abs\ju=
\ans{\ang{\xi_1,k_1},\dots,\ang{\xi_\nu,k_\nu}}$ and
$ 
\jta(\xi_i,k_i)=\ty{\xi_i}{m_i}{s_i}
$ 
for $i=1,\dots,\nu$, so that $\ju\leq\jta$. 

Consider the set $\cD$ of all \mus s 
$\Phi=\sid{\vpi_{\xi m}}{\xi<\jvt}{m<\om}\in \ms\bbp$ 
such that 
%$\Phi^J\cle\Phi$ and 
\ben
\nenu
%\atc\atc\atc\atc
\itla{1nen}
there is a number $\bah>h$ 
\imar{1nen}
and strings $\bas_i\in 2^{\bah}$ satisfying  
$s_i\su\bas_i$ and     
$\vys{\vpi_{\xi_im_i}}=\bah+1$ for $i=1,\dots,\nu$;

\itla{2nen}
there is a \mut\ $\jsg\in\mt\bbp$ 
%$\jsg\leq\jta$ 
which occurs in 
\imar{2nen}
$\Phi$ (Definition~\ref{muts}) and satisfies  
conditions \ref{3nen}, \ref{4nen} below; 

\itla{3nen}
$\jsg(\xi_i,k_i)=\vpi_{\xi_i m_i}(\bas_i)$ for $i=1,\dots,\nu$; 

\itla{4nen}
$\jsg$ directly forces $\rc\nin[T]$, where 
\imar{4nen}
$T=\bigcup_{s\in 2^{\bah}}\vpi_{\et M}(s)$.
\een

\ble
\lam{Ka}
$\cD$ is dense in\/ $\ms\bbp$. 
% above\/ $\Phi^J$. 
\ele
\bpf
By Lemma~\ref{xr}, 
it suffices to prove that for any
\mus\ 
$\Phi=\sid{\vpi_{\xi m}}{\xi<\jvt}{m<\om}\in\ms\bbp$ 
which already satisfies \ref{1nen} by means of a 
number $\bah$ and strings $\bas_i\in 2^{\bah}$, 
$1\le i\le\nu$,   
there is a \mus\ 
$\Phi'\in\cD$ which reduces $\Phi$.

Let 
$\xh=2^{\bah}$ (a number) and let 
$\ans{t_1,\dots,t_{\xh}}=2^{\bah}=\ens{t\in\bse}{\lh t=\bah}$. 
We suppose that the enumeration is chosen so that 
$t_i=\bas_i$ for $i=1,\dots,\mu$.
%be an enumeration of all 
%strings of length $H-1$, such that $t_i=\bas_i$ for each 
%$i\le\mu$. 
Let $\ell_i=k_i$ whenewer $1\le i\le \mu$.  
If $\mu+1\le n\le\xh$ then let 
$$
\ell_n=n+1+\tmax_{1\le i\le \nu}\ens{k_i}{\xi_i=\et}\,,
$$ 
so 
that pairs of the form $\ang{\et,\ell_n}$, $n\ge\mu+1$, 
do not belong to $\abs\jta$.  

Consider a \mut\ 
$\jro\in\mt\bbp$, defined so that 
\bit 
\item 
$\abs\jro=\abs\jta   
\cup\ens{\ang{\et,\ell_n}}{\mu+1\le n\le \xh}$; 

\item
$\jro(\xi_i,k_i)=\vpi_{\xi_i m_i}(\bas_i)
\text{ \ for all \ $i=1,\dots,\nu$}$;
%, so that $\jro\leq\jta$ since $s_i\sq\bas_i$ for all $i$; 

\item
$\jro(\et,\ell_n)= \vpi_{\et M}(t_n)$ for all 
$n$, $\mu+1\le n\le\xh$ --- note that by construction 
the equality $\jro(\et,\ell_i)= \vpi_{\et M}(t_i)$ 
also holds for $i=1,\dots,\mu$, being just a 
reformulation of $\jro(\xi_i,k_i)=\vpi_{\xi_i m_i}(\bas_i)$.
\eit
By the density of sets $D(k)$, there exists a \mut\  
$\jsg\in\mt\bbp\yt \jsg\leq\jro$, which 
directly forces $\rc\ne\rpi_{\et\ell_n}$ for all 
$n=1,\dots,\xh$ ---
including $\rc\ne\rpi_{\et k_i}$ for $i=1,\dots,\mu$. 
Then there are strings $u,v_1,\dots,v_{\xh}\in\bse$ 
such that $u$ is incompatible in $\bse$ with each 
$v_{n}$ and $\jsg$ directly forces each of the formulas 
\bce
$u\su\rc$, \ \ and \ \ 
$v_n\sq\rpi_{\et\ell_n}$   
for all $n$, $1\le n\le \xh$. \ 
%including $v_1\sq\rpi_{\et M}$. 
\ece
However $\jsg$ directly forces $v_n\sq\rpi_{\et\ell_n}$ iff 
$v_n\sq\roo{\jsg(\et,\ell_n)}$. 
We conclude that $\jsg$ directly forces\/ $\rc\nin[T]$,
where $T=\bigcup_{1\le n\le \xh}\jsg(\et,\ell_n)$.

Now let 
$\Phi'=\sid{\vpi'_{\xi m}}{\xi<\jvt}{m<\om}\in\ms\bbp$ 
be defined as follows. 
\ben
\Renu
\itla{step1}
we let 
\imar{step1}
$\vpi'_{\xi_i m_i}(\bas_i)=\jsg(\xi_i,k_i)$  \ 
for \ $i=1,\dots,\nu$;

\itla{step2}
if $\mu+1\le n\le \xh$ then let 
$\vpi'_{\et M}(t_n)=\jsg(\et,\ell_n)$ --- the equality is 
also true for $n\le\mu$ by \ref{step1};

\itla{step3}
if $\ang{\xi,m}\in\abs\Phi$, $s\in\bse,$ and 
$\lh s<\vys{\vpi_{\xi m}}$ 
(that is, $\vpi_{\xi m}(s)$ is defined), 
but $\vpi'_{\xi m}(s)$ is {\ubf not} defined by \ref{step1} 
and \ref{step2}\snos
{That is, except for the triples  
$\ang{\xi,m,s}=\ang{\xi_i,m_i,\bas_i}$ 
and $\ang{\et,M,t_n}$.}, 
then we keep 
$\vpi'_{\xi m}(s)=\vpi_{\xi m}(s)$;

\itla{step4}
for any $\ang{\xi,k}\in\abs\jsg\bez\abs\jro$ 
add to $\abs{\Phi'}$ 
a pair $\ang{\xi,m}\nin\abs{\Phi}$ and define 
$\vys{\vpi'_{\xi m}}=1$, $\vpi'_{\xi m}(\La)=\jsg(\xi,k)$ 
--- to make sure that $\jsg$ occurs in $\Phi'$.
\een
By construction, the \mus\ $\Phi'\in \ms\bbp$ 
%obtained this way 
reduces $\Phi$, the \mut\ $\jsg$ 
occurs in $\Phi'$ by \ref{step4} and satisfies 
$\jsg\leq\jro$. 
Finally to check \ref{4nen} note that by construction 
$\bigcup_{1\le n\le \xh}\jsg(\et,\ell_n)=
\bigcup_{s\in 2^{\bah}}\vpi'_{\et M}(s)$.    
Thus $\Phi'\in\cD$, as required.
\epF{Lemma}

Come back to the proof of the theorem.
It follows from the lemma that there is an index $j$ 
such that the system 
$\Phi^j=\sid{\vpt j\xi m}{\xi<\jvt}{m<\om}$ 
belongs to $\cD$. 
Let this be witnessed by a number $H>h$, a collection of 
strings $\bas_i\in\gh$ ($1\le i\le\nu$), and 
a \mut\ $\jsg\in\mt\bbp$,  
so that conditions 
\ref{1nen}, \ref{2nen}, \ref{3nen}, \ref{4nen} 
are satisfied for $\Phi=\Phi^j$ and $\jsg$. 
Then, for instance, 
$\vpt j{\xi_i}{m_i}(\bas_i)=\ty{\xi_i}{m_i}{\bas_i}$
(see Definition~\ref{dPhi}(iii)). 
However $\jsg(\xi_i,k_i)=\vpt j{\xi_i}{m_i}(\bas_i)$ by \ref{3nen} 
while $\jta(\xi_i,k_i)=\ty{\xi_i}{m_i}{s_i}$ by the construction, 
and $s_i\su\bas_i$. 
It follows that 
\vyk{
$\jsg(\et,M)\sq\jta(\et,M)$. 

Then, for instance, $\vpt j\et M(\bas)= \ty \et M \bas$
(see Definition~\ref{dPhi}(iii)). 
However $\jsg(\et,M)=\vpt j\et M(\bas)$ by \ref{3nen} 
while $\jta(\et,M)=\ty \et M s$ by the construction, 
and $s\su\bas$. 
It follows that $\jsg(\et,M)\sq\jta(\et,M)$. 

Similarly, 
$\jsg(\et,M')\sq\jta(\et,M')$ and 
$\jsg(\za,N)\sq\jta(\za,N)$, 
and hence in general we have
} 
$\jsg\leq\jta$.

Finally consider a \mut\ $\jv\in\mt\bbu$,  
defined so that $\abs\ju=\abs\jsg$, 
$
\ju(\xi_i,k_i)=\tf{\xi_i}{m_i}{\bas_i}
$
for $i=1,\dots,\nu$,  
and if 
$\ang{\xi,k}\in\abs\jsg \bez \ens
{\ang{\xi_i,k_i}}{1\le i\le \nu}$ 
then let $\jv(\xi,k)$ be any tree in $\dU_{\xi k}$ 
satisfying $\jv(\xi,k)\sq\jsg(\xi,k)$ 
(we refer to Lemma~\ref{uu2}).

Recall that by construction $s_i\su\bas_i$ for all $i$.  
It follows that $\jv\leq\ju$. 
On the other hand, $\jv\leq\jsg$, therefore 
$\jv$ directly forces $\rc\nin[T]$ by \ref{4nen}, where 
$T=\bigcup_{s\in 2^{H}}\vpt j\et M(s)
=\bigcup_{s\in 2^{H}}\ty\et Ms$.
And finally by definition  
$U=\tfi\et M\sq\bigcup_{s\in 2^{H}}\ty\et Ms$, so 
$\jv$ directly forces $\rc\nin[U]$, as required. 
\epf 

%***************
%}

\vyk{
\bte
\lam{L}
In the assumptions of Definition~\ref{dPhi}, suppose that\/ 
a set\/ $D\in\cM$, $D\sq\mt\bbp$, is pre-dense in\/ 
$\mt\bbp$. 
Then it remains pre-dense in\/ $\bbpu$.
\ete
}

\vyk{
\bpf
By construction $U\sq \tfi\et M$ for some $M<\om$; 
thus we can assume that simply $U=\tfi\et M$. 
The indices 
$\et\yi M$ are fixed in the proof. 
We can assume by Lemma~\ref{uu2} that $\ju\in\mt\bbu$. 
The support 
$\abs\ju=\ans{\ang{\xi_1,k_1},\dots,\ang{\xi_\nu,k_\nu}}
\sq\jvt\ti\om$
is a finite set ($\nu<\om$), 
and if $i=1,\dots,\nu$ then, as $\ju\in\mt\bbu$, 
there is a string $s_i$ and a number $m_i$ such that 
$\ju(\xi_i,k_i)=\tf{\xi_i}{m_i}{s_i}$.
By obvious reasons including Lemma~\ref{disj}, 
we can assume that 
\ben
\aenu
\itla{ass2}
if $i\ne i'$ and $\xi_i=\xi_{i'}$ then $k_i\ne k_i'$;
\imar{ass2}

\itla{ass3}
$s_i\ne s_{i'}$ whenever $i\ne i'$, 
\imar{ass3}
and there is $h<\om$ 
such that $\lh{s_i}=h\yd\kaz i$;

\itla{ass1}
let $H=2^h$, then 
$\nu\ge H$, 
\imar{ass1}
$\xi_1=\dots=\xi_{H}=\et$, $m_1=\dots=m_{H}=M$, and 
$\ans{s_1,\dots,s_H}=\ens{s\in\bse}{\lh s=h}$;

\itla{ass4}
if $H<i\le\nu$ then 
\imar{ass4}
$\ang{\xi_i,m_i}\ne \ang{\et,M}$.
\een
Comment to \ref{ass3}:
if $s_i\su s'_i\in\bse$ for all $i$, and $\ju'\in\mt\bbu$ 
is defined so that $\abs{\ju'}=\abs{\ju}$ and 
$\ju'(\xi_i,k_i)=\tf{\xi_i}{m_i}{s'_i}$ for all $i=1,...,\nu$, 
then $\ju'\leq \ju$.
Therefore, if we prove the theorem for $\ju'$ then it implies 
the result for $\ju$ as well. 

In these assumptions, by Lemma~\ref{uu4} 
there exists a \mut\  
$\jv\in\mt\bbu\yt \jv\leq\ju$, which 
directly forces $\rc\ne\rpi_{\et k_i}$ for all 
$i=1,\dots,H$; recall that $\et=\xi_i$ whenever 
$1\le i\le H$.  
Then there are strings $u,t_1,\dots,t_{H}\in\bse$ 
such that $u$ is incompatible in $\bse$ with each 
$t_{i}$ and $\jv$ directly forces each of the formulas 
\bce
$u\su\rc$, \ and \  
$t_n\sq\rpi_{\et k_i}$   
for $1\le i\le H$. 
\ece
However $\jv$ directly forces $t_i\sq\rpi_{\et k_i}$ iff 
$t_i\sq\roo{\jv(\et,k_i)}$. 
We conclude that $\jv$ directly forces\/ $\rc\nin[T]$,
where $T=\bigcup_{1\le i\le H}\jv(\et,k_i)$.

Now let 
$\Phi'=\sid{\vpi'_{\xi m}}{\xi<\jvt}{m<\om}\in\ms\bbp$ 
be defined as follows. 
\ben
\Renu
\itla{step1}
we let 
\imar{step1}
$\vpi'_{\xi_i m_i}(\bas_i)=\jsg(\xi_i,k_i)$  \ 
for \ $i=1,\dots,\nu$;

\itla{step2}
if $\mu+1\le n\le \gh$ then let 
$\vpi'_{\et M}(t_n)=\jsg(\et,\ell_n)$ --- the equality is 
also true for $n\le\mu$ by \ref{step1};

\itla{step3}
if $\ang{\xi,m}\in\abs\Phi$, $s\in\bse,$ and 
$\lh s<\vys{\vpi_{\xi m}}$ 
(that is, $\vpi_{\xi m}(s)$ is defined), 
but $\vpi'_{\xi m}(s)$ is {\ubf not} defined by \ref{step1} 
and \ref{step2}\snos
{That is, except for the triples  
$\ang{\xi,m,s}=\ang{\xi_i,m_i,\bas_i}$ 
and $\ang{\et,M,t_n}$.}, 
then we keep 
$\vpi'_{\xi m}(s)=\vpi_{\xi m}(s)$;

\itla{step4}
for any $\ang{\xi,k}\in\abs\jsg\bez\abs\jro$ 
add to $\abs{\Phi'}$ 
a pair $\ang{\xi,m}\nin\abs{\Phi}$ and define 
$\vys{\vpi'_{\xi m}}=1$, $\vpi'_{\xi m}(\La)=\jsg(\xi,k)$ 
--- to make sure that $\jsg$ occurs in $\Phi'$.
\een
By construction, the \mus\ $\Phi'\in \ms\bbp$ obtained 
this way satisfies $\Phi\cle\Phi'$, the \mut\ $\jsg$ 
occurs in $\Phi'$ by \ref{step4} and satisfies 
$\jsg\leq\jro$. 
Finally to check \ref{4nen} note that by construction 
$\bigcup_{1\le n\le \gh}\jsg(\et,\ell_n)=
\bigcup_{s\in 2^{\bah-1}}\vpi'_{\et M}(s)$.    
Thus $\Phi'\in\cD$, as required.

%\epF{Lemma}

Come back to the proof of the theorem.
It follows from the lemma that there is an index $j$ 
such that the system 
$\Phi^j=\sid{\vpt j\xi m}{\xi<\jvt}{m<\om}$ 
belongs to $\cD$. 
Let this be witnessed by a number $H>h$, a collection of 
strings $\bas_i\in\gh$ ($1\le i\le\nu$), and 
a \mut\ $\jsg\in\mt\bbp$,  
so that conditions 
\ref{1nen}, \ref{2nen}, \ref{3nen}, \ref{4nen} 
are satisfied for $\Phi=\Phi^j$ and $\jsg$. 
Then, for instance, 
$\vpt j{\xi_i}{m_i}(\bas_i)=\ty{\xi_i}{m_i}{\bas_i}$
(see Definition~\ref{dPhi}(iii)). 
However $\jsg(\xi_i,k_i)=\vpt j{\xi_i}{m_i}(\bas_i)$ by \ref{3nen} 
while $\jta(\xi_i,k_i)=\ty{\xi_i}{m_i}{s_i}$ by the construction, 
and $s_i\su\bas_i$. 
It follows that 
\vyk{
$\jsg(\et,M)\sq\jta(\et,M)$. 

Then, for instance, $\vpt j\et M(\bas)= \ty \et M \bas$
(see Definition~\ref{dPhi}(iii)). 
However $\jsg(\et,M)=\vpt j\et M(\bas)$ by \ref{3nen} 
while $\jta(\et,M)=\ty \et M s$ by the construction, 
and $s\su\bas$. 
It follows that $\jsg(\et,M)\sq\jta(\et,M)$. 

Similarly, 
$\jsg(\et,M')\sq\jta(\et,M')$ and 
$\jsg(\za,N)\sq\jta(\za,N)$, 
and hence in general we have
} 
$\jsg\leq\jta$.

Finally consider a \mut\ $\jv\in\mt\bbu$,  
defined so that $\abs\ju=\abs\jsg$, 
$
\ju(\xi_i,k_i)=\tf{\xi_i}{m_i}{\bas_i}
$
for $i=1,\dots,\nu$,  
and if 
$\ang{\xi,k}\in\abs\jsg \bez \ens
{\ang{\xi_i,k_i}}{1\le i\le \nu}$ 
then let $\jv(\xi,k)$ be any tree in $\dU_{\xi k}$ 
satisfying $\jv(\xi,k)\sq\jsg(\xi,k)$ 
(we refer to Lemma~\ref{uu2}).

Recall that by construction $s_i\su\bas_i$ for all $i$.  
It follows that $\jv\leq\ju$. 
On the other hand, $\jv\leq\jsg$, therefore 
$\jv$ directly forces $\rc\nin[T]$ by \ref{4nen}, where 
$T=\bigcup_{s\in 2^{H-1}}\vpt j\et M(s)
=\bigcup_{s\in 2^{H-1}}\ty\et Ms$.
And finally by definition  
$U=\tfi\et M\sq\bigcup_{s\in 2^{H-1}}\ty\et Ms$, so 
$\jv$ directly forces $\rc\nin[U]$, as required. 
\epf 

}

\parf{The product forcing} 
\las{jfor}

In this section, 
{\ubf we argue in $\rL$, the constructible universe.}
Let $\lel$ be the canonical wellordering of $\rL$.

\bdf
[in $\rL$]
\lam{uxi}
We define, by induction on $\al<\omi$, sequences   
$\bbu^\al=\sis{\dU^\al_\xi}{\xi<\al}$,  
$\bbp^\al=\sis{\dP^\al_\xi}{\xi<\al}$ 
of countable sets of trees 
$\dU^\al_\xi\yi \dP^\al_\xi$ in $\ptf$, 
as follows.

First of all, we let $\dP^\al_\al=0$ and $\dU^\al_\al=\dP_0$ 
(see Example~\ref{p0}) 
for all $\al$; 
note that the terms $\dP^\al_\al\yi\dU^\al_\al$ 
do not participate in 
the sequences $\bbp^\al$ and $\bbu^\al$.

{\ubf The case \boldmath$\al=0$.} 
Let $\bbp^0=\bbu^0=\La$ 
(the empty sequence).

{\ubf The step.} 
Suppose that $0<\la<\omi$, and  
$\bbu^\al$, $\bbp^\al$ as above
are already defined for every $\al<\la$. 
%, so that $(\dag)$ holds. 
Let $\cM_\la$ be the least model $\cM$ of $\zfc'$ of the form 
$\rL_\ka\yt\ka<\omi$, containing 
$\sis{\bbu^\al}{\al<\la}$ 
and $\sis{\bbp^\al}{\al<\la}$, 
and such that 
$\la<\omim$  
and $\dU^\al_\xi\yi \dP^\al_\xi$ are countable in $\cM$ 
for all $\xi<\al<\la$.

We first define a sequence 
$\bbp^\la=\sis{\dP^\la_\xi}{\xi<\la}$ so that 
$\dP^\la_\xi=\bigcup_{\xi\le\al<\la}\dU^\al_\xi$
for all $\xi<\la$.
In particular if $\la=\al+1$ then      
$\dP^{\al+1}_\xi=\dP^\al_\xi\cup \dU^\al_\xi$ for all 
$\xi<\al+1$ 
(because $\dP^\al_\xi=\bigcup_{\xi\le\al'<\al}\dU^{\al'}_\xi$ 
at the previous step),  
and, for $\xi=\al$, 
$\dP^{\al+1}_{\al}=\dP^\al_\al\cup \dU^\al_\al=\dP_0$ 
(see above). 
Thus $\bbp^{\al+1}$ is the
extension of $\bbp^\al\lor\bbu^\al$ 
(see Section~\ref{mul}) 
by the default assignment $\dP^{\al+1}_{\al}=\dP_0$.
For instance, $\bbp^1=\ang{\dP_0}$.

Thus a sequence $\bbp^\la=\sis{\dP^\la_\xi}{\xi<\la}$ 
is defined. 

To define $\bbu^\la$ and 
accomplish the step, 
let $\dphi=\sis{\Phi^j}{j<\om}$ 
be the \dd\lel least sequence of 
\mus s $\Phi^j\in\ms{\bbp^\la}$, 
\dd\cle increasing and generic over $\cM_\la$, 
and let $\bbu^\la=\sis{\dU^\la_\xi}{\xi<\la}$ 
be defined, on the base of this sequence, as 
in Definition~\ref{dPhi}. 

After the sequences 
$\bbu^\al=\sis{\dU^\al_\xi}{\xi<\al}$ and 
$\bbp^\al=\sis{\dP^\al_\xi}{\xi<\al}$, 
and the model $\cM_\al$, have been 
defined for all $\al<\omi$, we let 
$\dP_\xi=\bigcup_{\xi\le\al<\omi}\dU^\al_\xi$ 
for all $\xi<\omi$, and let 
$\bbp=\bbp^{\omi}=\sis{\dP_\xi}{\xi<\omi}$.  
The set $\mt\bbp$ of all \dd\bbp multitrees 
(Definition~\ref{muts}) 
will be our principal forcing notion. 
\edf

\vyk{
Note that $\vt:\omi\to\omi$ is a strictly increasing function, 
but $\ran\vt$ can be a proper subset of $\omi$. 
Nevertheless we define a quasi-inverse 
$\ivt:\omi\onto\omi$ 
such that $\ivt(\xi)=\al$ whenever $\vt(\al)\le\xi<\vt(\al+1)$.
(The underlining 
indicates that this is not a true inverse map.)
Thus $\xi<\vt(\al)\eqv \ivt(\xi)<\al$, and the sets of trees 
$\dU^\al_\xi\yi\dP^\al_\xi$ are defined by Definition~\ref{uxi} 
whenever $\ivt(\xi)<\al<\omi$.
}

\bpro
%[in $\rL$]
\lam{uxip}
The sequences\/ $\sis{\bbu^\al}{\al<\omi}$, 
$\sis{\bbp^\al}{\al<\omi}$ belong to\/ 
$\id{\hc}1$.\qed
\epro

\bre
\lam{a<g}
If $\al<\ga\le\omi$ then the sets $\mt{\bbp^\al}$ 
and $\mt{\bbp^\ga}$ of \mut s are formally disjoint. 
However we can naturally embed the former in the latter. 
Indeed each \mut\ 
$\jta=\sid{T_{\xi k}}{\xi<\al}{k<\om}\in\mt{\bbp^\al}$
can be identified as an element of $\mt{\bbp^\ga}$ by the 
default extension $T_{\xi k}=\bse$ whenever 
$\al\le\xi<\ga$. 
With such an identification, we can assume that 
$\mt{\bbp^\al}\sq\mt{\bbp^\ga}\sq\mt\bbp$, and similarly 
$\mt{\bbp^\la}=\bigcup_{\al<\la}\mt{\bbu^\al}$ for all 
limit $\la$, and the like.
\ere

\ble
%[in $\rL$]
\lam{jden}
If\/ $\al<\omi$ and a set\/ $D\in\cM_\al\yt D\sq \mt{\bbp^\al}$ 
is pre-dense in\/ $\mt{\bbp^\al}$ then it remains pre-dense in\/ 
$\mt{\bbp}$. 

Therefore\/ 
$\mt{\bbu^\al}$ is pre-dense in\/ $\mt{\bbp}$.
\ele
\bpf
By induction on $\ga\yd \xi\le\ga<\omi$, 
if $D$ is pre-dense in $\mt{\bbp^\ga}$ then it 
remains pre-dense in 
$\mt{\bbp^\ga\lor \bbu^\ga}$ 
by Lemma~\ref{uu4}, hence in $\mt{\bbp^{\ga+1}}$ too by 
constructions. 
Limit steps including the step $\omi$ are obvious. 

To prove the second part, note that 
$\mt{\bbu^\al}$ is dense in $\mt{\bbp^\al\lor \bbu^\al}$ 
by Lemma~\ref{uu2}, therefore, pre-dense in 
$\mt{\bbp^{\al+1}}$,  
and $\mt{\bbu^\al}\in\cM_{\al+1}$.
\epf

\bcor
%[in $\rL$]
\lam{xiden}
If\/ $\xi<\al<\omi$ then the set\/ $\dU^\al_\xi$ 
is pre-dense in\/ $\dP_\xi$.
\ecor
\bpf
Let $T\in \dP_\xi$. 
Consider a \mut\ $\jta\in\mt{\bbp}$ defined so that 
$\jta(\xi,0)=T$ and $\jta(\eta,k)=\bse$ whenever 
$\ang{\et,k}\ne\ang{\xi,0}$. 
By Lemma~\ref{jden} $\jta$ is compatible in $\mt\bbp$ 
with some $\ju\in\mt{\bbu^\al}$. 
We conclude that $T$ is compatible in $\dP_\xi$ with 
$U=\ju(\xi,0)\in \dU^\al_\xi$.
\epf

\ble
%[in $\rL$]
\lam{club}
If\/ $X\sq\hc=\rL_{\omi}$ then the set\/ $W_X$ of all 
ordinals\/ $\al<\omi$ such that\/ 
$\stk{\rL_\al}{X\cap\rL_\al}$ is an elementary submodel of\/  
$\stk{\rL_{\omi}}{X}$ and\/ $X\cap\rL_\al\in\cM_\al$ 
is unbounded in\/ $\omi$.
More generally, if\/ $X_n\sq\hc$ for all\/ $n$ 
then the set\/ $W$ of all 
ordinals\/ $\al<\omi$, such that\/ 
$\stk{\rL_\al}{\sis{X_n\cap\rL_\al}{n<\om}}$ 
is an elementary submodel of\/  
$\stk{\rL_{\omi}}{\sis{X_n}{n<\om}}$ 
and\/ $\sis{X_n\cap\rL_\al}{n<\om}\in\cM_\al$, 
is unbounded in\/ $\omi$.
\ele
\vyk{
\bpf
Let $\xi_0<\omi$. 
By standard arguments, there are ordinals $\xi<\la<\omi$, 
$\xi>\xi_0$, such 
that $\stk{\rL_\la}{\rL_\xi,X\cap\rL_\xi}$ is an elementary 
submodel of $\stk{\rL_{\om_2}}{\rL_{\omi},X}$.
Then $\stk{\rL_\xi}{X\cap\rL_\xi}$ is an elementary submodel of   
$\stk{\rL_{\omi}}{X}$, of course. 
Moreover, $\xi$ is uncountable in $\rL_\la$, hence 
$\rL_\la\sq\cM_\xi$. 
It follows that $X\cap\rL_\xi\in\cM_\xi$ since 
$X\cap\rL_\xi\in\rL_\la$ by construction.
The second claim does not differ much.
\epf
}

\bpf
Let $\al_0<\omi$. 
Let $M$ be a countable elementary submodel of $\rL_{\om_2}$ 
containing $\al_0\yi\omi\yi X$, 
and such that $M\cap\hc$ is transitive. 
Let $\phi:M\onto\rL_\la$ be the Mostowski collapse, and 
let $\al=\phi(\omi)$. 
Then $\al_0<\al<\la<\omi$ and $\phi(X)=X\cap\rL_\al$ by the 
choice of $M$. 
It follows that $\stk{\rL_\al}{X\cap\rL_\al}$ is an elementary 
submodel of $\stk{\rL_{\omi}}{X}$.
Moreover, $\al$ is uncountable in $\rL_\la$, hence 
$\rL_\la\sq\cM_\al$. 
We conclude that $X\cap\rL_\al\in\cM_\al$ since 
$X\cap\rL_\al\in\rL_\la$ by construction.

The second claim does not differ much.
\epf

\bcor
%[in $\rL$, = Lemma 6 in \cite{jenmin}]
\lam{ccc}
The forcing\/ $\mt\bbp$ satisfies CCC.
\ecor
\bpf
Suppose that $A\sq\mt\bbp$ is a maximal antichain. 
By Lemma~\ref{club}, there is an ordinal $\al$ such that 
$A'=A\cap\mt{\bbp^\al}$ is a maximal antichain in 
$\mt{\bbp^\al}$ and $A'\in\cM_\al$. 
But then $A'$ remains pre-dense, therefore, maximal, in the 
whole set $\mt\bbp$ by Lemma~\ref{jden}. 
It follows that $A=A'$ is countable.
\epf

\parf{The extension: non-uniformizable set and Theorem~\ref{Tun}} 
\las{mod}

Working in terms of Definition~\ref{uxi}, 
we consider the set $\mt\bbp\in\rL$     
as a forcing notion over $\rL$. 
It is equal to the finite-support product 
$\prod_{\xi<\omi}{\dP_\xi}\lom$, 
which also can be understood as 
the finite-support product 
$\prod_{\xi<\omi,\, k<\om}{\dP_{\xi k}}$, where each 
$\dP_{\xi k}$ is equal to one and the same 
$\dP_\xi=\bigcup_{\xi\le\al<\omi}\dU^\al_\xi$ 
of Definition~\ref{uxi}.  

We make use of this forcing to {\ubf prove Theorem~\ref{Tun}}.
 
\ble
[= Lemma 7 in \cite{jenmin}]
\lam{mod1}
Let\/ $\xi<\omil$.
A real\/ $x\in\dn$ is\/ $\dP_\xi$-generic over\/ $\rL$ iff\/ 
$x\in Z_\xi=\bigcap_{\xi<\al<\omil}\bigcup_{U\in\dU^\al_\xi}[U]$. 
\ele
\bpf
All sets $\dU^\al_\xi$ are pre-dense in $\dP_\xi$ by 
Corollary~\ref{xiden}. 
On the other hand, if $A\sq\dP_\xi$, $A\in\rL$ is a maximal 
antichain in $\dP_\xi$, then easily $A\sq\dP^\al_\xi$ for some 
$\al$, $\xi<\al<\omil$, by Corollary~\ref{ccc}.
But then every tree $U\in\dU^\al_\xi$ satisfies 
$U\sqf\bigcup A$ by Lemma~\ref{uu3}, so that 
$\bigcup_{U\in \dU^\al_\xi}[U]\sq \bigcup_{T\in A}[T]$.
\epf       

\bcor
%[= Corollary 9 in \cite{jenmin}]
\lam{mod2}
In any generic extension of\/ $\rL$ with the same\/ $\omi$, 
the set 
$$
P=\ens{\ang{\xi,x}}{\xi<\omil\land x\in\dn \,
\text{is \dd{\dP_\xi}generic over $\rL$}}\sq\omil\ti\dn
$$
is\/ $\ip\hc1$, and\/ $\ip12$ in terms of a usual 
coding system of ordinals\/ $<\omi$ by reals. 
\ecor
\bpf
Use Lemma~\ref{mod1} and Proposition~\ref{uxip}.
\epf  

\bdf
\lam{gg}
From now on, let $G\sq\plo$ be a set \dd{\mt\bbp}generic over 
$\rL$.
Note that $\omi^{\rL[G]}=\omil$ by Corollary~\ref{ccc}.

If $\xi<\omil$ and $k<\om$ then let 
$G_{\xi k}=\ens{\jta(\xi,k)}{\jta\in G}$, 
so that each $G_{\xi k}$ is 
\dd{\dP_\xi}generic over $\rL$ and 
$X_{\xi k}=\bigcap_{T\in G_{\xi k}}[T]$ is 
a singleton $X_{\xi k}=\ans{x_{\xi k}}$ whose only element 
$x_{\xi k}\in\dn$ is 
a real \dd{\dP_\xi}generic over $\rL$.
\edf

The whole extension $\rL[G]$ is then equal to 
$\rL[\sis{x_{\xi k}}{\xi<\omil,\, k<\om}]$, 
and our goal is now to prove that 
it contains no \dd{\dP_\xi}generic reals except for 
the reals $x_{\xi k}$.

\ble
[in the assumptions of Definition~\ref{gg}] 
\lam{only}
If\/ $\xi<\omil$ and\/ $x\in\rL[G]\cap\dn$ then\/ 
$x\in\ens{x_{\xi k}}{k<\om}$ iff\/ 
$x$ is a\/ \dd{\dP_\xi}generic real over\/ $\rL$.
\ele
\bpf
Otherwise there is a \mut\ $\jta\in\mt\bbp$ and 
a \dd{\mt\bbp}real name 
$\rc=\sis{\kc ni}{n<\om,\,i=0,1}\in\rL$ 
such that $\jta$ \dd{\mt\bbp}forces that $\rc$ is 
\dd{\dP_\xi}generic over $\rL$ 
while $\mt\bbp$ forces $\rc\ne\rpi_{\xi k}$, $\kaz k$. 
(Recall that $\rpi_{\xi k}$ is a \dd{\mt\bbp}name for 
$x_{\xi k}$.)

The sets $C_n=\kc n0\cup\kc n1$ are pre-dense in 
$\mt\bbp$. 
It follows from Lemma~\ref{club} that there is an ordinal 
$\la$, $\xi<\la<\omi$, such that 
%$\jta\in{\dP_\la}\lom$, 
each set 
$C'_n=C_n\cap \mt{\bbp^\la}$ is pre-dense in $\mt{\bbp^\la}$, 
and the sequence $\sis{\xc ni}{n<\om,\,i=0,1}$ belongs to 
$\cM_\la$, where $\xc ni=C'_n\cap \kc ni$ --- 
then $C'_n$ is pre-dense in $\mt{\bbp}$, too, 
by Lemma~\ref{jden}. 
Thus we can assume that in fact $C_n=C'_n$, that is, 
$\rc\in\cM_\la$ and $\rc$ is a \dd{\mt{\bbp^\la}}name.

Further, as $\mt\bbp$ forces that $\rc\ne\rpi_{\xi k}$, 
the set $D_{k}$ 
of all \mut s $\jsg\in\mt\bbp$ which directly force 
$\rc\ne\rpi_{\xi k}$, is dense in $\mt\bbp$ --- for every $k$. 
Therefore, still by Lemma \ref{club}, 
% and \ref{jden}, 
we may assume that the same ordinal $\la$ as above satisfies 
the following: 
each set $D'_{k}=D_{k}\cap\mt{\bbp^\la}$ 
is dense in $\mt{\bbp^\la}$. 
%and belongs .$ 
%(And pre-dense in $\plo$ by Lemma~\ref{jden}).

Applying Theorem~\ref{K} with 
$\bbp=\bbp^\la$, 
$\bbu=\bbu^\la$, 
$\jvt=\la$, $\et=\xi$,
%and $\dP\cup\dU=\dP_{\la+1}$, 
we conclude that for each 
$U\in\dU_\xi^\la$ the set $Q_U$ of all \mut s 
$\jv\in \mt{\bbu^\la}$ which directly force $\rc\nin[U]$, 
is dense in $\mt{\bbu^\la\lor\bbp^\la}$, therefore, 
pre-dense in $\mt{\bbp^{\la+1}}$.
As obviously $Q_U\in\cM_{\la+1}$, we further conclude that 
$Q_U$ is pre-dense in $\mt\bbp$ 
by Lemma~\ref{jden}.
Therefore $\mt\bbp$ forces 
$\rc\nin\bigcup_{U\in \dU_\xi^\la}[U]$, 
hence, forces that $\rc$ is not \dd{\dP_\xi}generic, 
by Lemma~\ref{mod1}.
But this contradicts to the choice of $\jta$.
\epf

\ble
[in the assumptions of Definition~\ref{gg}] 
\lam{sym}
If\/ $\xi<\omil$ and\/ $k<\om$ then\/ 
\ben
\renu
\itla{sym1}
$x_{\xi k}\nin 
\rL[\sis{x_{\et\ell}}{\ang{\et,\ell}\ne\ang{\xi,k}}]$, 

\itla{sym2}
$x_{\xi k}$ is not\/ 
$\od(\sis{x_{\et\ell}}{\et\ne\xi,\,k<\om})$ in\/ 
$\rL[G]$.\qed
\een
\ele
\bpf
\ref{sym1} is a usual property of product forcing, 
while to prove \ref{sym2} we need to make use of the 
fact that by construction the 
\dd\xi part of the forcing is itself a finite-support 
product of countably many copies of $\dP_\xi$.
\epf

\bpri
[non-uniformizable $\ip\hc1$ set]
\lam{nonu1}
Arguing in the assumptions of Definition~\ref{gg},
we consider, 
in $\rL[G]=\rL[\sis{x_{\xi k}}{\xi<\omil,\,k<\om}]$, 
the set $P$ of Corollary~\ref{mod2}.
First of all $P$ is $\ip\hc1$ in $\rL[G]$ by Corollary~\ref {mod2}.
Further, it follows from Lemma~\ref{only} that 
$$
P=\ens{\ang{\xi,x_{\xi k}}}{\xi<\omil\land k<\om}\,,
$$
and hence all vertical cross-sections of $P$ are countable. 
And by Lemma~\ref{sym} it is not ROD uniformizable since 
any real in $\rL[G]$ belongs to a submodel of the form 
$\rL[\sis{x_{\xi k}}{\xi<\za,\, k<\om}]$, where $\za<\omil$.
\epri

\bpri
[non-uniformizable $\ip12$ set]
\lam{nonu2}
To get a non-uniformizable $\ip12$ set in $\dn\ti\dn$ on 
the base of the abovedefined set $P\sq\omil\ti\dn,$ we 
make use of a usual coding of countable ordinals by 
reals. 
Let $\wo\sq\dn$ be the $\ip11$ set of codes, 
and for $w\in\wo$ 
let $\abs w<\omi$ be the ordinal coded by $w$.
We consider
$$
P'=\ens{\ang{w,x}\in\wo\ti\dn}{\ang{\abs w,x}\in P}\,, 
$$
a $\ip12$ set in $\rL[G]$.
Suppose towards the contrary that, in $\rL[G]$, $P'$ 
is uniformizable by a ROD set $Q'\sq P'$.
As $\omil=\omi$ by Corollary~\ref{ccc}, for any $\xi<\omi$ 
there is a code $w\in\wo\cap\rL$ with $\abs w=\xi$. 
Let $w_\xi$ be the \dd\lel least of those. 
Then
$$
Q=\ens{\ang{\xi,x}\in P}{\ang{w_\xi,x}\in Q'}  
$$
is a ROD subset of $P$ which uniformizes $P$, contrary to 
Example~\ref{nonu1}.
\epri

\qeD{Theorem~\ref{Tun}}

\parf{Non-separation model} 
\las{ns}

Here we prove Theorem~\ref{Tsep}. 
The model we use will be defined on the base of the model 
$\rL[G]=\rL[\sis{x_{\xi k}}{\xi<\omil,\, k<\om}]$
of Definition~\ref{gg}, of the form 
$\gN_\Xi=\rL[\sis{x_{\xi 0}}{\xi\in\Xi}]$, 
where $\Xi\sq\omil$ will be a generic subset of $\omil$, 
so that, strictly speaking, $\gN_\Xi$ is not going to be a 
submodel of $\rL[G]$. 

To define $\Xi$, we recall first of all that the ordinal 
product $2\nu$ is considered as the ordered sum of $\nu$ 
copies of $2=\ans{0,1}$. 
Thus if $\nu=\la+m$, where $\la$ is a limit ordinal or $0$ 
and $m<\om$, then $2\nu=\la+2m$ and $2\nu+1=\la+2m+1$. 

Now let $\dQ=3^{\omil}$ with finite support, so that a typical 
element of $\dQ$ is a partial map $q:\omil\to 3={0,1,2}$ with 
a finite domain $\dom q\sq\omil$; 
this is a version of the Cohen forcing, of course.  

\bdf
[in the assumptions of Definition~\ref{gg}] 
\lam{ggh}
Let $H\sq\dQ$ be a set generic over $\rL[G]$. 
It naturally yiels a Cohen-generic map $F_H:\omil\to3$. 
Let 
$$
\bay{ccccrrrcccccccc}
A_H&=&\ens{\nu<\omil}{F_H(\nu)=0}\,,&&  
B_H&=&\ens{\nu<\omil}{F_H(\nu)=1}\,,\\[1ex]
D_H&=&\ens{\nu<\omil}{F_H(\nu)=2}\,,&&\text{and}
\eay
$$  
$$
\Xi_H=\ens{2\nu}{\nu\in A_H\cup D_H}\cup
\ens{2\nu+1}{\nu\in B_H\cup D_H}\,.
$$
We consider the model 
$\gN_H=\rL[\sis{x_{\xi 0}}{\xi\in\Xi_H}]$.
Let $\hc(H)=(\hc)^{\gN_H}$.
\edf

Note that $\gN_H$ is not a submodel of $\rL[G]$ since the 
set $\Xi_H$ does not belong to $\rL[G]$; but 
$\gN_H\sq\rL[G][H]$, of course.

\bte
[in the assumptions of Definition~\ref{ggh}]
\lam{Tsep'}
It is true in\/ $\gN_H$ that\/ $A_H$ and \/ $B_H$ are 
disjoint\/ $\ip{\hc(H)}2$ sets
not 
separable by disjoint\/ $\fs\hc2$ sets.
\ete

\bpri
[non-separable $\ip13$ sets]
\lam{nons}
%Before we begin the proof of this theorem, let us show 
%how it implies Theorem~\ref{Tsep}.
In the notation of Example~\ref{nonu2}, let 
$$
X=\ens{w_\xi}{\xi\in A_H}
\quad\text{and}\quad
Y=\ens{w_\xi}{\xi\in B_H}\,.
$$
The sets $X,Y\sq\wo\cap\rL$ are $\ip{\hc(H)}2$   
together with $A_H$ and $B_H$, and hence $\ip13$, 
and $X\cap Y=\pu$. 
Suppose towards the contrary that $X',Y'\sq\dn$ are 
disjoint sets in $\fs13$, hence in $\fs{\hc(H)}2$, 
such that $X\sq X'$ and $Y\sq Y'$. 
Then 
$$
A=\ens{\xi<\omil}{w_\xi\in X'}
\quad\text{and}\quad
B=\ens{\xi<\omil}{w_\xi\in Y'}
$$
are disjoint sets in $\fs{\hc(H)}2$, 
and we have $A_H\sq A$ and $B_H\sq B$ by construction, 
contrary to Theorem~\ref{Tsep'}.
\epri
 
The proof of Theorem~\ref{Tsep'} involves the following 
result which will be established in the next section.
Theorem \ref{gex} esentially says that the coding structure 
in $\rL[G]$ described in Section~\ref{mod} survives a further 
Cohen-generic extension.

\bte
[Cohen-generic stability]
\lam{gex}
In the assumptions of Definition~\ref{ggh}$:$  
\ben
\renu
\itla{gex1}
if\/ $\xi<\omil$ and\/ $x\in\rL[G][H]\cap\dn$ then\/ 
$x\in\ens{x_{\xi k}}{k<\om}$ iff\/ 
$x$ is a\/ \dd{\dP_\xi}generic real over\/ $\rL\;;$

\itla{gex2}
if\/ $\xi<\omil$ and\/ $k<\om$ then\/ 
$x_{\xi k}\nin 
\rL[\sis{x_{\et\ell}}{\ang{\et,\ell}\ne\ang{\xi,k}}][H]\;;$

\itla{gex3}
if\/ $\xi<\omil$ and\/ $k<\om$ then\/ 
$x_{\xi k}$ is not\/ 
$\od(\sis{x_{\et\ell}}{\et\ne\xi,\,k<\om},H)$ in\/ 
$\rL[G][H]$.
\een
\ete 

\bpf[Theorem~\ref{Tsep'} modulo Theorem \ref{gex}] 
That $A_H\cap B_H=\pu$ is clear. 
To see that, say, $A_H$ is $\ip{\hc(H)}2$ in $\gN_H$, 
prove that the equality
$$
A_H=\ens{\nu<\omi}{\neg\:\sus x\,P(2\nu+1,x)}
$$
holds in $\gN_H$, where $P$ is the $\ip\hc1$ set of 
Corollary~\ref{mod2}.
(For $B_H$ it would be $P(2\nu,x)$ in the displayed formula.)

First suppose that $\nu<\omil$, $\xi=2\nu+1$, 
$x\in\gN_H\cap\bn$, and $P(\xi,x)$ holds 
in $\gN_H$; prove that $\nu\nin A_H$. 
By definition $x$ is \dd{\dP_{\xi}}generic over $\rL$. 
Then $x=x_{\xi k}$ for some $k$ by 
Theorem~\ref{gex}\ref{gex1}. 
Therefore $k=0$ and $\xi$ has to belong to $\Xi_H$ 
by Theorem~\ref{gex}\ref{gex2}. 
But then $\nu\in B_H\cup D_H$, so $\nu\nin A_H$, as required. 

To prove the converse, suppose that $\nu\nin A_H$, so that 
$\nu\in B_H\cup D_H$. 
Then $\xi=2\nu+1\in\Xi_H$, and hence $x=x_{\xi0}\in\gN_H$.
We conclude that $\ang{\xi,x}=\ang{2\nu+1,x}\in P$ 
by Lemma~\ref{only}, as required.

Finally, {\ubf to prove the non-separability}, 
suppose towards the contrary that, in $\gN_H$, $A_H$ and 
$B_H$ are separable by a pair of disjoint $\fs\hc2$ sets 
$A,B\sq\omi=\omil$.
These sets are defined in the set 
$\hc(H)=(\hc)^{\gN_H}$ by $\ip{}2$ 
formulas, resp., $\vpi(a,\xi)\yi\psi(b,\xi)$, with real 
parameters $a,b\in\gN_H\cap\dn.$ 
Let $\la<\omil$ be a limit ordinal such that 
$a,b\in\rL[\sis{x_{\xi 0}}{\xi\in\Xi_H\cap\la}]$, and let 
$\sg,\tau\in\rL[G]$ be \dd\dQ real names such that 
$a=\sg[H]$ and $b=\tau[H]$, which depend on 
$\sis{x_{\xi 0}}{\xi\in\Xi_H\cap\la}$ only.

If $K\sq\dQ$ is a set \dd\dQ  generic over $\rL[G]$ 
(\eg, $K=H$), then let 
$$
\ahh{K}=\ens{\xi<\omil}
{\vpi(\sg[{K}],\xi)^{\hc({K})}}\,,\;\,
\bhh{K}=\ens{\xi<\omil}
{\psi(\tau[{K}],\xi)^{\hc({K})}}\,,
$$
so that by definition $A_H\sq A=\ahh H$, $B_H\sq B=\bhh H$, 
and $\ahh H\cap\bhh H=\pu$.
Fix a condition $q_0\in H$ which forces, over $\rL[G]$,  
that
$A_\uH\sq \ahh\uH$,  
$B_\uH\sq \bhh\uH$, and 
$\ahh \uH\cap\bhh \uH=\pu$,
where $\uH$ is the canonical name for $H$.
We may assume that $\dom{q_0}\sq\la $ as well, 
for otherwise just increase $\la$.

Now let $\xi_0$ be any ordinal with $\la\le\xi_0<\omi$. 
Consider three sets $H_0\yi H_1\yi H_2\sq\dQ$, generic 
over $\rL[G]$ and containing $q_0$, 
whose generic maps $F_{H_i}:\omil\to3$ 
satisfy $F_{H_i}(\xi_0)=i$ and 
$F_{H_0}(\xi)=F_{H_1}(\xi)=F_{H_2}(\xi)$ 
for all $\xi\ne\xi_0$.

Then $\sg[H_0]=\sg[H_2]$, $\tau[H_0]=\tau[H_2]$,
and $\Xi_{H_2}=\Xi_{H_0}\cup\ans{2\xi_0+1}$, hence,  
$\gN_{H_0}\sq\gN_{H_2}$. 
It follows by Shoenfield that 
$\ahh{H_0}\sq\ahh{H_2}$ and   
$\bhh{H_0}\sq\bhh{H_2}$, hence
$$
A_{H_2}\sq A_{H_0}\sq\ahh{H_0}\sq\ahh{H_2}\,,\;\,   
B_{H_2}=B_{H_0}\sq\bhh{H_0}\sq\bhh{H_2}\,,\;\,
\ahh{H_2}\cap\bhh{H_2}=\pu
$$
by the choice of $q_0$. 
We conclude that $\xi_0\in \ahh{H_2}$, just because 
$\xi_0\in A_{H_0}$ by the choice of $H_0$.
And we have $\xi_0\in \bhh{H_2}$ by similar reasons. 
Thus $\ahh{H_2}\cap\bhh{H_2}\ne\pu$, 
contrary to the above.
The contradiction ends the proof.\vtm

\epF{Theorems~\ref{Tsep'} and~\ref{Tsep} modulo 
Theorem \ref{gex}}

\parf{The proof of the Cohen-generic stability theorem} 
\las{ns}

Here {\ubf we prove Theorem \ref{gex}.}
We concentrate on Claim \ref{gex1} of the theorem since 
claims \ref{gex2}, \ref{gex3} are established by the same 
routine product-forcing arguments outlined in 
the proof of Lemma~\ref{sym}.

First of all, let us somewhat simplify the task. 
It is known that every real in a \dd\dQ generic extension 
belongs to a simple \dd\bse generic extension 
(that is, a Cohen-generic one) of the same model.
That is, it suffices to prove this:

\ble
[in the assumptions of Definition~\ref{gg}] 
\lam{only+}
If\/ $a\in\dn$ is \dd\bse generic over $\rL[G]$, $\xi<\omil$, 
and a real $x\in\rL[G][a]\cap\dn$ 
is \dd{\dP_\xi}generic over $\rL[G]$ 
then $x=x_{\xi k}$ for some $k$.
\ele
\bpf
Coming back to Definition~\ref{uxi}, we conclude that 
the sequence $\dphi$ there is generic not only over 
$\cM_\la$ but also over $\cM_\la[a]$ by the product forcing 
theorems.
It follows that Lemma~\ref{jden} also is true in $\rL[a]$ 
for all sets $D\in\cM_\al[a]$, and so are Lemma~\ref{club} 
(for models $\rL_{\omi}[a]$ and $\rL_\al[a]$) 
and Corollaries \ref{xiden} and \ref{ccc}.
This enables us to prove Lemma~\ref{only} for all reals 
$x\in\rL[G][a]$, and we are done.
\epf

\qeD{Theorem \ref{gex}}

\vyk{
The proof of the lemma (compare to Lemma~\ref{only}) 
is based on the following 
Theorem~\ref{K+}, a strengthening of Theorem~\ref{K}. 
We let $\mtp\bbp=\mt\bbp\ti\bse$, with the componentwise 
order $\ang{\jta,q}\leq\ang{\jta',q'}$ iff $\jta\leq\jta'$ 
and $q\sq q'$; $\mtp\bbp$ adds a sequence of the form 
$\ang{\sis{x_{\xi k}}{\xi<\jvt,\,k<\om},a}$, where 
$\sis{x_{\xi k}}{\xi<\jvt,\,k<\om}$ is \dd{\mt\bbp}generic 
while $a\in\dn$ is \dd\bse generic.
%To formulate it, we need to accordingly modify the notation 
%involved in Theorem~\ref{K}.

We'll use $q,r$ in this section to denote   
strings in $\bse.$

We come back to Section~\ref{mul}. 
Let $\vt\le\omi$ and let $\bbp=\sis{\dP_\xi}{\xi<\vt}$ 
be a sequence of sets $\dP_\xi\in\ptf$.
If $\ang{\jta,q}\in\mtp\bbp$ then we define 
$[\jta,q]=[\jta]\ti [q]\sq2^{\vt\ti\om}\ti\dn$, where 
$[q]=\ens{x\in\dn}{q\su x}$ while $[\jta]\sq2^{\vt\ti\om}$ 
was defined in Section~\ref{mul}.

Return to Section~\ref{saway}. 
A \rit{\dd{\mtp\bbp}real name} is a system 
$\rc=\sis{\kc ni}{n<\om,\, i<2}$ of sets $\kc ni\sq\mtp\bbp$ 
such that each set $C_n=\kc n0\cup \kc n1$ is 
pre-dense in $\mtp\bbp$ 
and if $\jsg\in \kc n0$ and $\jta\in \kc n1$ then $\jsg,\jta$ are 
incompatible in $\mtp\bbp$.

If $\xi<\vt$ and $k<\om$ then 
a \dd{\mtp\bbp}real name 
$\rpi_{\xi k}=\sis{\jcp ni\xi k}{n<\om,\,i<2}$
is defined so that each set $\jcp ni\xi k$ contains a 
single condition  
$\ang{\jr ni\xi k,\La}\in\mtp\bbp$, where  
$\jr ni\xi k\in\mt\bbp$ was introduced by Definition~\ref{proj} 
and $\La\in\dQ$ is the empty string.
Then $\rpi_{\xi k}$ is a \dd{\mtp\bbp}real name of the real 
$x_{\xi k}$, the $(\xi,k)$th 
term of a \dd{\mtp\bbp}generic sequence 
$\ang{\sis{x_{\xi k}}{\xi<\jvt,\,k<\om},a}$. 

The notions of a condition $\ang{\jta,q}\in\mtp\bbp$ 
directly forcing some elementary formulas are introduced 
as in Section~\ref{saway}, with obvious corrections.

\bte
\lam{K+}
In the assumptions of Definition~\ref{dPhi}, suppose that\/ 
$\et<\vt$, $\rc= 
\sis{C_m^i}{m<\om,\,i<2}\in\cM$ is a\/ 
\dd{\mtp\bbp}real name, 
and for every\/ $k$ the set
$$
D(k)=\ens{\ang{\jta,q}\in\mtp\bbp}
{\ang{\jta,q}\text{ directly forces }
\rc\ne\rpi_{\et k}}
$$
is dense in\/ $\mtp\bbp$. 
Let\/ $\ang{\ju,q}\in\mtp\bbpu$ and\/ 
$U\in \dU_\et$.
Then there is\/ a condition\/ 
$\ang{\jv,r}\in\mtp\bbu\yd \ang{\jv,r}\leq\ang{\ju,q}$, 
which directly forces\/ $\rc\nin[U]$.
\ete

\bpf[Lemma~\ref{only+}]
Similar to Lemma~\ref{only+}, on the base of 
Theorem~\ref{K+}.
\epf
}

\bibliographystyle{plain}
{\small
%\bibliography{u,su}
%\input{su.bbl}

%
}

\end{document}